\documentclass[10pt]{amsart}
\newcommand*\colorgraphics{}
\usepackage[a4paper]{geometry}

\usepackage{ifthen}
\usepackage{pifont}
\usepackage{url}
\usepackage{graphicx}
\usepackage[utf8]{inputenc}
\newtheorem{prop}[subsection]{Proposition}
\newtheorem{cor}[subsection]{Corollary}
\newtheorem{thm}[subsection]{Theorem}
\newtheorem{lemma}[subsection]{Lemma}
\newtheorem{defi}[subsection]{Definition}

\numberwithin{equation}{section}

\DeclareMathOperator{\PNG}{PNG}
\DeclareMathOperator{\RSK}{RSK}
\DeclareMathOperator{\GUE}{GUE}
\DeclareMathOperator{\Aztec}{A}
\DeclareMathOperator{\Loz}{L}
\DeclareMathOperator{\Herm}{GUE}
\DeclareMathOperator{\Ge}{Ge}
\DeclareMathOperator{\Var}{Var}
\DeclareMathOperator{\Vol}{Vol}
\DeclareMathOperator{\Real}{Re}

\begin{document}
\newboolean{dum}
\setboolean{dum}{false}
\newcommand{\comment}[1]{\ifthenelse{\boolean{dum}}{
{\par\noindent\Huge\ding{46}} \fbox{#1}\par}{}}
\newcommand{\gcn}{\ensuremath{C_{\mathbf{GC,\mathbb{N}}}}}
\newcommand{\gc}{\ensuremath{C_{\mathbf{GC}}}}

\title{Eigenvalues of GUE minors}
\author{Kurt Johansson}
\address{Institutionen för Matematik, 
Swedish Royal Institute of Technology (KTH), 100 44 Stockholm, Sweden}
\email{kurtj@math.kth.se}

\author{Eric Nordenstam}
\address{Institutionen för Matematik, 
Swedish Royal Institute of Technology (KTH), 100 44 Stockholm, Sweden}
\email{eno@math.kth.se}

\bibliographystyle{alpha}
\comment{
$ $Id: thesisart1.tex,v 1.4 2009/01/08 12:43:38 enord Exp $ $
}
\thanks{This work was supported by the Göran Gustafsson foundation (KVA)}

\begin{abstract}
Consider an infinite random matrix $H=(h_{ij})_{0<i,j}$
picked from the Gaussian Unitary Ensemble (GUE). 
Denote
its main minors by $H_i=(h_{rs})_{1\leq r,s\leq i}$
and let the $j$:th largest eigenvalue of 
$H_i$ be $\mu^i_j$. 
We show that 
the configuration of all these eigenvalues $(i,\mu_j^i)$
form a determinantal point process on 
$\mathbb{N}\times\mathbb{R}$.

Furthermore we show that this process can be obtained as the scaling limit 
in random tilings of the Aztec diamond close to the boundary.
We also discuss the corresponding limit for random 
lozenge tilings of a hexagon. 
\end{abstract}
\maketitle

This version of this article differs from the one published
in Electronic Journal of Probability in that the errors listed
in the separate erratum have been corrected. 

\section{Introduction}
The distribution of eigenvalues induced by some measure on matrices 
has been the study of random matrix theory for decades. 
These distributions have been found to be universal in the
sense that they turn up in various unrelated problems,
some of which do not contain a matrix in any obvious way,
or contain a matrix that does not look like a random matrix.
In this article, we propose to study the eigenvalues of the minors
of a random matrix, and argue that this distribution also is 
universal in some sense by showing that it is the scaling limit 
of three apparently unrelated discrete models. 

The largest eigenvalues of minors of GUE-matrices have been studied in 
\cite{baryshnikov:gue}, connecting these to a certain queueing model.
It is a special case of the very general class
of models analysed in \cite{johansson:detproc}.
The large $N$ limit of this model will yield 
the distribution of all the eigenvalues of a GUE-matrix and its minors.

This process will turn out also to be the
scaling limit of a point process related to random tilings of 
the Aztec diamond, studied in \cite{johansson:aztec}
and of a process related to random lozenge-tilings of a hexagon,
studied in \cite{johansson:rhombus}.

\subsection{Eigenvalues of the GUE}
Consider the following point process  on $\Lambda=\mathbb{N}\times\mathbb{R}$.
There is  a point at $(n, \mu)$ iff the $n$:th main minor of $H$, i.e. $H_n$, 
has an eigenvalue $\mu$. 
We will call this process the \emph{GUE minor process\/.}
A central result in this article is that this process
is a determinantal point process with a certain kernel 
$K^{\Herm}$.

For details of what it means for a point process to be 
determinantal, see section~\ref{sec:determ-point-proc}.
An explicit expression for this kernel is given in the next definition.

\begin{defi}
\label{def:herm}
The  \emph{GUE minor kernel} is
\begin{equation*}
K^{\Herm}(r,\xi; \,s,\eta)=-\phi(r,\xi;\,s,\eta) +
\sum_{j=-\infty}^{-1}\sqrt{\frac{(s+j)!}{(r+j)!}}h_{r+j}(\xi)h_{s+j}(\eta) e^{-(\xi^2+\eta^2)/2}, 
\end{equation*}
where $\phi(r,\xi;\,s,\eta)=0$ when $r\leq s$ and 
\begin{multline*}
\phi(r,\xi;\,s,\eta)=\frac{(\xi-\eta)^{r-s-1}\sqrt{2^{r-s}}}{(r-s-1)!} e^{\frac{1}{2} (\eta^2-\xi^2)}H(\xi-\eta)\\
-\frac{\frac 12 e^{(\eta^2-\xi^2)/2}}{\sqrt[4]{\pi}}
\sum_{j=-r}^{-(s+1)}
\frac{h_{r+j}(\xi)\sqrt{2^{-s-j}}}{\sqrt{(r+j)!}(-s-j-1)!}
\int_{\eta}^\infty (t-\eta)^{-s-j-1}e^{-t^2}dt
\end{multline*}
for $r>s$.
\end{defi}
Here,  $h_k(x)=2^{-k/2}(k!)^{-1/2}\pi^{-1/4} H_k(x)$ are the 
Hermite polynomials normalised
so that \[\int h_i(x)h_j(x) e^{-x^2} dx=\delta_{ij},\]
$h_k\equiv 0$ when $k<0$ and $H$ is the Heaviside function defined by 
\begin{equation}
\label{eq:heaviside}
H(t)=\begin{cases} 
1 & \text{for $t\geq 0$} \\
0 & \text{for $t<0$.}
\end{cases}
\end{equation}

\begin{thm}
\label{thm:gue}
The GUE minor process is determinantal
with  kernel $K^{\Herm}$.
\end{thm}
This will be proved in section~\ref{sec:gue-minor-kernel}. 
\subsection{Aztec Diamond}
\label{sec:aztec-diamond}
\begin{figure}[p]
\includegraphics[width=12.5cm]{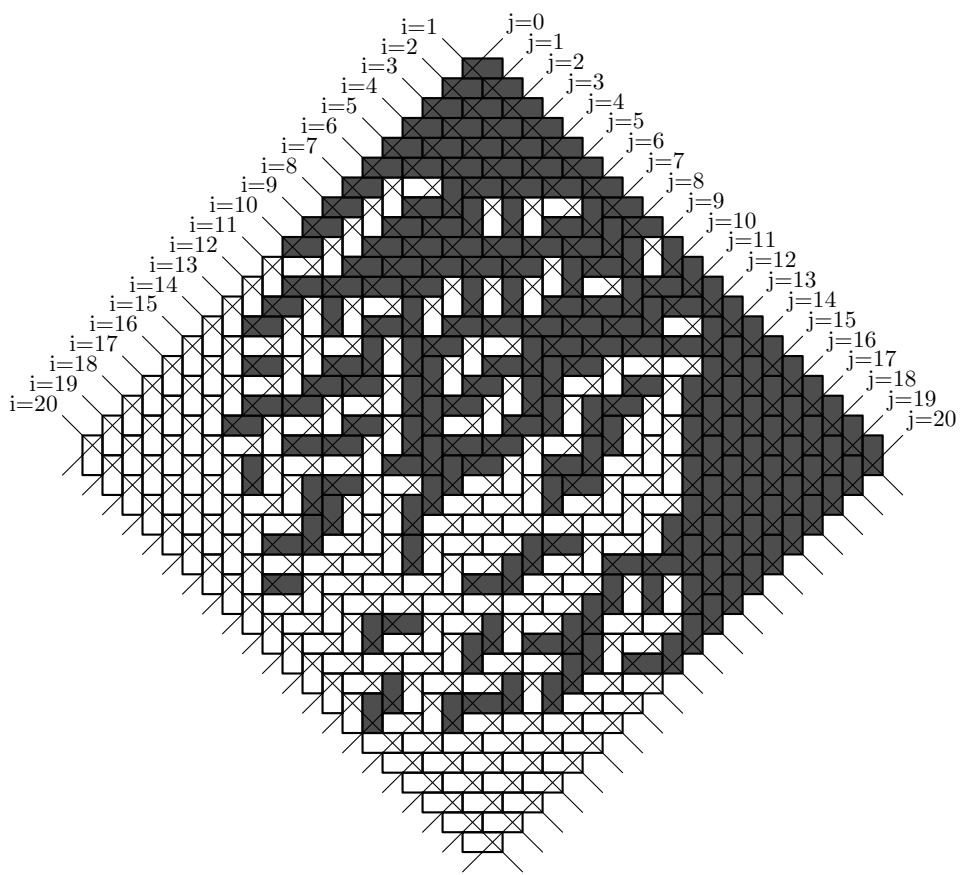}
\caption{An Aztec Diamond of size 20 with N and E type dominoes shaded.}
\label{fig:aztec}
\end{figure}
The Aztec diamond of size $N$  is the largest region of the plane that
is the union of squares with corners in lattice points and
is contained in the region $|x|+|y|\leq N+1$, see figure \ref{fig:aztec}.

It can be covered with $2\times1$ dominoes in
$2^{N(N+1)/2}$ ways,~\cite{elkies:alternating_i,elkies:alternating_ii}.
Probability distributions on the set of  all these
possible tilings have been studied in several 
references, for example \cite{johansson:aztec,propp:aztec}. 
Typical samples are characterized by having so called
frozen regions in the north, south, east and west, regions
where the tiles are layed out like brickwork. 
In the middle there is a disordered region, the so called
temperate region.
It is for example known that for large $N$,
the shape of the temperate region tends to a circle, 
see \cite{propp:arctic_circle} for precise statement.

The key to analyzing this model is to colour all squares black
or white in a checkerboard fashion. 
Let us chose colour white for the left
square on the top row. A horizontal tile is of type N, or north,
whenever its left square is white. All other horizontal tiles are
of type S, or south. 
Likewise, a vertical domino is of type W, or west,  precisely if its top square
is white. Other vertical dominoes are of type E. 

In figure \ref{fig:aztec}, tiles of type  N and E  have been shaded.
Notice that along the  line $i=1$, there is precisely
one white tile, and its position is a stochastic variable that
we denote $x_1^1$. 
Along the line $i=2$ there are precisely two
white tiles, at positions $x_1^2$ and $x_2^2$ respectively, etc.
In general, let $x^i_k$ denote the $j$-coordinate of the $k$:th 
white tile along line $i$. 
These white points can be considered a particle configuration,
and this particle configuration uniquely determines the tiling.
It is shown in \cite{johansson:aztec} that this process is a
determinantal point process on 
$\mathbf{N}^2=\{1,2,\dots,N\}^2$, and
the kernel is computed.

We show that this particle process, properly rescaled,
converges weakly to the distribution for eigenvalues of GUE described
above. More precisely we have the following theorem that will
be proved in section~\ref{sec:aztec-diamond-1}.
\begin{thm}
\label{thm:main_aztec_thm}
Let $\mu^i_j$ be the eigenvalues of a GUE matrix and its minors.
For each $N$, let $\{x_j^i\}$ be the position of the particles, as
defined above, 
in a random tiling of the Aztec Diamond of size $N$.
Then for each continuous function of compact 
support $\phi:\mathbb{N}\times\mathbb{R}\rightarrow\mathbb{R}$,
with $0\leq \phi\leq 1$,
\begin{equation*}
\mathbb{E}\left[
\prod_{i,j} (1-\phi(i,\mu^i_j))\right] =
\lim_{N\rightarrow\infty} 
\mathbb{E}\left[
\prod_{i,j}(1-\phi(i,\frac{x^i_j-N/2}{\sqrt{N/2}}
))\right].
\end{equation*}
\end{thm}

\subsection{Rhombus Tilings}

\begin{figure}[t]
\includegraphics{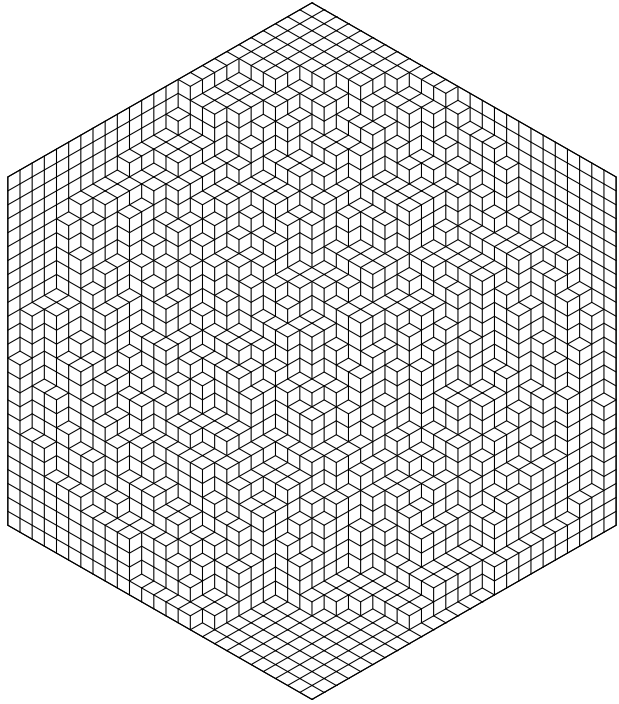}
\caption{Lozenge tiling of a hexagon.}
\label{fig:large_hexagon}
\end{figure}
 
\begin{figure}[p]
\colorgraphics
\includegraphics[width=12.5cm]{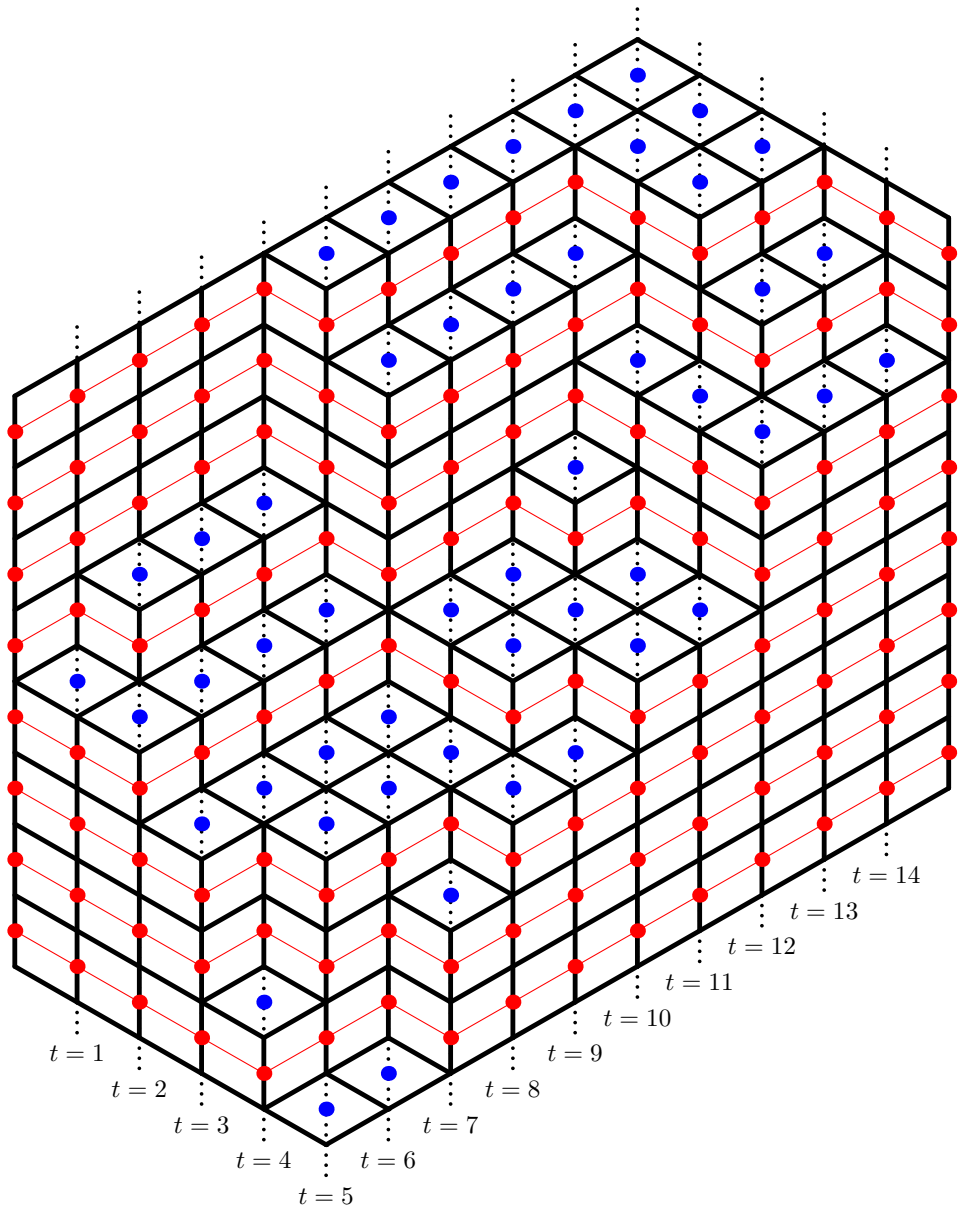}
\caption{Tiled  hexagon with sides $a=8$, $b=5$ and $c=10$.
 The so called horizontal rhombuses are marked with a blue dot. }
\label{fig:small_hexagon}
\end{figure}

Consider an $(a,b,c)$-hexagon, i.e. a hexagon  with 
side lengths $a$, $b$, $c$, $a$, $b$, $c$.  
It can be covered by rhombus-shaped tiles with
angles $\pi/3$ and $2\pi/3$ and
side length $1$, so called \emph{lozenges}. 
The number of possible such tilings is
\begin{equation*}
\prod_{i=1}^a\prod_{j=1}^b\prod_{k=1}^c
\frac{i+j+k-1}{i+j+k-2}.
\end{equation*}
This formula was proved by Percy MacMahon (1854-1928), 
see \cite[page 401]{stanley:ec2} for historical remarks.

Thus, we can  chose a tiling randomly, 
each possible tiling assigned equal  probability.
A typical such tiling is shown in figure \ref{fig:large_hexagon}.
Just like in the case of the Aztec diamond, 
there are frozen regions in the corners of the shape and
a disordered region in the middle. 
It has been shown, that when $a=b=c=N\rightarrow \infty$,
this so called temperate region,
tends to a circle, see \cite{cohn:shape_boxed_partition} for
precise statement and other similar results.

Equivalently, consider $a$
 simple, symmetric, random walks, started at 
positions $(0, 2j)$, $1\leq j\leq a$. 
At each step
in discrete time, each walker moves up or down, 
with equal probability. 
They are  conditioned never to intersect 
and to end at positions $(c+b, c-b+2j)$. 
Figure \ref{fig:small_hexagon}, the red lines 
illustrate such a family of walkers, and shows the correspondence
between this process and tilings of the hexagon.
These red dots in the figure  define a point process.
 \cite{johansson:rhombus} shows that uniform measure on
tilings of the hexagon (or equivalently, uniform measure on 
the  possible configurations of simple, symmetric, random walks) 
induces a measure on this point process that is determinantal, and
computes the kernel.

We will show, in theorem~\ref{thm:blue-dots},
that the complement of this process, the blue dots
in the figure, is also a determinantal process and compute its
kernel. 

Let us introduce some notation. 
Observe that along the line $t=1$, there is exactly one blue dot.
Let its position be $x_1^1$. 
Along line $t=2$ there are two blue dots, at positions $x_1^2$ 
and $x_2^2$ respectively, and so on. All these $x_j^i$ are
stochastic variables, and they are of course not independent 
of each other. 

We expect that  the scaling limit of the process
$\{x_j^i\}_{i,j}$, 
as the size of the hexagon tends to infinity, is the GUE minor 
process with kernel $K^{\Herm}$.   
More precisely, 
let $\mu_j^i$ be the eigenvalues of a GUE matrix and its minors. 
Then
for each continuous function of compact support 
$\phi:\mathbb{N}\times \mathbb{R}\rightarrow\mathbb{R}$, $0\leq \phi\leq 1$, 
\begin{equation}
\label{eq:15}
\mathbb{E}\left[\prod_{i,j}(1-\phi(i,\mu^i_j))\right]=
\lim_{N\rightarrow\infty}
\mathbb{E}\left[\prod_{i,j}(1-\phi(i,\frac{x^i_j-N/2}{\sqrt{N/2}}))\right].
\end{equation}
We will outline a proof of this result by going to the limit in the formula 
for the correlation kernel, which involves the Hahn polynomials.
A complete proof requires some further estimates of these polynomials. 

The GUE minor process has also been obtained as a 
limit at ``turning points'' 
in a 3D partition model by 
Okounkov and Reshetikhin~\cite{okounkov:birth}. 
We expect that the GUE minor process should be the universal limit in 
random tilings where the disordered region touches the boundary.

Acknowledgement: We thank A.~Okounkov for helpful comments and for 
sending the preprint~\cite{okounkov:birth}. 

\section{Determinantal point processes}
\label{sec:determ-point-proc}
Let $\Lambda$ be a complete separable metric space with some
reference measure $\lambda$. For example $\mathbb{R}$ with 
Lebesgue measure or $\mathbb{N}$ with counting measure. 
Let $M(\Lambda)$
be the space  of integer-valued and locally finite measures on $\Lambda$. 
A \emph{point process} $x$ is a probability measure  on $M(\Lambda)$.
For example, let $x$ be a point process. 
A realisation  $x(\omega)$ is an element of  $M(\Lambda)$.  It will
assign positive measure to certain points, 
$\{x_i(\omega)\}_{1\leq i\leq N(\omega)}$, sometimes called \emph{particles},
or just points in the process. 
In the processes that we will study the number of particles
in a compact set will have a uniform upper bound.

Many point processes can be specified by giving their 
\emph{correlation functions}, 
$\rho_n:\Lambda^n\rightarrow\mathbb{R}$, $n=1,\dots,\infty$.
We will not go into the precise definition of these
 or when a process is uniquely determined
by its correlation functions. For that we refer to any 
or all of the following references:
\cite[Ch. 9.1, A2.1]{daley:pointproc}, \cite{johansson:randmatdetproc,soshnikov:detrandpointfield}.

Suffice it to say that correlation functions have the 
following useful property.
For any bounded measurable function $\phi$
with bounded support $B$, satisfying
\begin{equation}
\label{eq:0}
\sum_{n=1}^\infty \frac{||\phi||_\infty^n}{n!}
\int_{B^n} \rho_n(x_1,\dots,x_n) \; d^nx < \infty
\end{equation} 
the following holds:
\begin{equation}
\label{eq:1}
\mathbb{E}[\prod_i(1+\phi(x_i(\omega)))]=
\sum_{n=0}^\infty \frac{1}{n!}
\int_{\Lambda^n}
\phi(x_1)\cdots \phi(x_n)
\rho_n(x_1,\dots,x_n) \,d^n\lambda.
\end{equation}
Correlation functions are thus useful in computing various 
expectations. 
For example, if $A$ is some set and $\chi_A$ is
the characteristic function of that set, then 
$1-\mathbb{E}[\prod(1-\chi_A(x_i))]$ is the probability
of at least one particle in the set $A$. If the correlation 
functions of a process exist and are known, this probability
can then readily be computed with the above formula.

We will study point processes of a certain type, namely those 
whose correlation functions exist and are of the form
\begin{equation*}
\rho_n(x_1,\dots,x_n) = \det[K(x_i;\,x_j)]_{1\leq i,j\leq n},
\end{equation*}
i.e. the $n$:th correlation  function is given
as a $n\times n$ determinant where $K:\Lambda^2\rightarrow \mathbb{R}$ 
is some, not necessarily smooth, measurable function. 
Such a process is called a
\emph{determinantal point process}\/ and the function $K$ is called
the \emph{kernel} of the point process.

Let $x^1$, $x^2$, \dots, $x^N$,\dots
be a sequence of point processes on $\Lambda$.
Say that $x^N$ assigns positive measure to the points 
$\{x^N_i(\omega)\}_{1\leq i\leq N^N(\omega)}$.
Then we say that this sequence of point processes \emph{converges weakly}
to a point process $x$, written $x^N\rightarrow x$, $N\rightarrow \infty$,
 if for any continuous function $\phi$ of compact support, $0\leq \phi\leq 1$,

\begin{equation}
\lim_{N\rightarrow\infty} 
\mathbb{E}\left[\prod_{i=1}^{N^N(\omega)}(1-\phi(x_i^N(\omega)))\right]
= \mathbb{E}\left[\prod_{i=1}^{N(\omega)}(1-\phi(x_i(\omega)))\right].
\end{equation}
The next proposition gives sufficient 
conditions for weak convergence of a sequence of 
determinantal processes 
in terms of the kernels. 
\begin{prop}
\label{thm:determ-point-proc}
Let  $x^1$, $x^2$, \dots, $x^N$,\dots
be a sequence of determinantal point processes, 
and let $x^N$ have correlation kernel $K^N$ satisfying 
\begin{enumerate}
\item
$K^N\rightarrow K$, $N\rightarrow\infty$ pointwise, for some function $K$,
\item the $K^N$ are 
 uniformly bounded on compact sets in $\Lambda^2$ and
\item 
For $C$ compact, there exists some number $n=n(C)$ such that
\begin{equation*} 
\det[K^N(x_i,x_j)]_{1\leq i,j\leq m} = 0
\end{equation*} if $m\geq n$.
\end{enumerate}
Then there exists some determinantal point process $x$ with correlation
functions $K$ such that  $x^N\rightarrow x$ weakly.
\end{prop}
\begin{proof}We start by showing that there exists
such a determinantal point process $x$.
In \cite{soshnikov:detrandpointfield}, the following necessary and 
sufficient conditions for the existence of a random point process 
with given correlation functions is given.
\begin{enumerate}
\item Symmetry.
\begin{equation*}
\rho_k(x_{\sigma(1)},\dots,x_{\sigma(k)})=\rho_k(x_1,\dots,x_k)
\end{equation*}
\item Positivity.
For any finite set of measurable bounded functions 
$\phi_k:\Lambda^k\rightarrow\mathbb{R}$,
$k=0,\dots,M$, with compact support, such that 
\begin{equation}
\phi_0 + \sum_{k=1}^M \sum_{i_1\neq \dots\neq i_k}\phi(x_{i_1},\dots,x_{i_k})
\geq 0
\end{equation}
for all $(x_1,\dots,x_M)\in I^M$ it holds that
\begin{equation}
\phi_0+\sum_{k=1}^N \int_{I^k} \phi_k(x_1,\dots,x_k) \rho_k(x_1,\dots,x_k)
\;dx_1\dots dx_n\geq 0.
\end{equation}
\end{enumerate}

The first condition is satisfied for all correlation functions coming
from determinantal kernels because permuting the rows and the columns
of a matrix with the same permutation leaves the determinant unchanged.
For the positivity condition consider the kernels $K^N$. They are
 kernels  of determinantal processes so 
\begin{equation}
\phi_0+\sum_{k=1}^M \int_{I^k} \phi_k(x_1,\dots,x_k) 
\det[K^N(x_i,x_j)]_{1\leq i,j\leq k}
\;dx_1\dots dx_n\geq 0.
\end{equation}
As $N\rightarrow\infty$, this converges to the same expression with 
$K$ instead of $K^N$ by Lebesgue's bounded convergence theorem 
with assumption (2). Positivity  of this expression for all $N$ 
then implies positivity of the limit.

So now we know that $x$ exists.  We need to show that $x^N\rightarrow x$
Take some test function $\phi:\Lambda\rightarrow \mathbb{R}$
with bounded support $B$.  
For this function we check the condition in~(\ref{eq:0}). 
The assumption (3) in this theorem implies that the
sum is a finite one. Also, $||\phi||_\infty\leq 1$. 
Assumption (2) is that the functions $K^N$ are uniformly bounded, so
in particular they are bounded  on $B^2$, so $\rho_k$ is bounded
on $B^k$.
The integral of a bounded function over a bounded
set is finite, so this is the finite sum of finite real numbers,
which is finite. 

Therefore, for each $N$, by (\ref{eq:1}), 
\begin{align}
\lim_{N\rightarrow\infty} \mathbb{E}&
\left[\prod_{i}(1-\phi(x_i^N(\omega)))\right]=\\
&=
\lim_{N\rightarrow\infty}
\sum_{n=0}^\infty
\frac{(-1)^n}{n!}
\int_{\Lambda^n} \prod_{i=1}^n \phi(x_i)
\det[K^N(x_i,x_j)]_{1\leq i,j\leq n}\, d^n\lambda(x) .
\intertext{Condition (3) guarantees that the sum is finite.
Lebesgue's bounded convergence theorem applies because
the support of $\phi$ is compact and the correlation functions are
bounded on compact sets.  Thus the limit exists and is}
&=\sum_{n=0}^\infty
\frac{(-1)^n}{n!}
\int_{\Lambda^n} \prod_{i=1}^n \phi(x_i)
\det[K(x_i,x_j)]_{1\leq i,j\leq n}\, d^n\lambda(x) \\
&=\mathbb{E}\left[\prod_i(1-\phi(x_i(\omega)))\right].
\end{align}
This implies that indeed $x^N\rightarrow x$, weakly, as $N\rightarrow\infty$.
\end{proof}

\section{The GUE Minor Kernel}
\label{sec:gue-minor-kernel}
\subsection{Performance Table}
Consider the following model. 
Let $\{w(i,j)\}_{(i,j)\in\mathbb{Z}_+^2}$, 
be independent geometric random variables with parameter $q^2$. 
I.e. there is one i.i.d. variable sitting at each integer 
lattice point in the first quadrant of the plane.
Let 
\begin{equation}
G(M,N) = \max_{\pi} \sum_{(i,j)\in\pi} w(i,j)
\end{equation}
where the maximum is over all up/right paths from $(1,1)$ to $(M,N)$.
The array $[G(M,N)]_{M,N\in\mathbb{N}}$ is called the 
\emph{performance table}.

Each such up/right path must pass through precisely one 
of $(M-1,N)$ and $(M,N-1)$, so it is true that
$G(M,N)=\max(G(M-1,N),G(M,N-1))+w(M,N)$.

It is known from \cite{baryshnikov:gue}, that 
$(G(N,1), G(N,2),\dots,G(N,M))$ for fixed $M$, properly rescaled,
jointly tends
to the distribution of $(\mu_1^1,\mu_1^2,\dots,\mu_1^M)$
as $N\rightarrow\infty$ in the sense of weak convergence of
probability measures.
We will show that it is possible to define 
stochastic variables in terms of the values $w$ 
that jointly converge weakly to the distribution of
all the eigenvalues $\mu_j^i$ of GUE-matrices.

\subsection{Notation}
We will use the following notation from \cite{baryshnikov:gue}.
\begin{enumerate}
\item $W_{M,N}$ is  set of $M\times N$ integer matrices.
\item $W_{M,N,k}$ is  set of $M\times N$ integer matrices whose entries sum
up to $k$.
\item $V_M=\mathbb{R}^{M(M-1)/2}$
where the components of each element $x$ are indexed in the following way.
\begin{equation*}
x=\begin{matrix}
x_1^1\\
\vdots & \ddots\\
x_1^{M-1} & \dots & x_{M-1}^{M-1}\\
x_1^{M} & \dots & x_{M}^{M-1} &  x_{M}^{M}\\
\end{matrix}
\end{equation*}
\item $\gc\subset V_M$ is the subset such that $x^i_{j-1}\geq x_{j-1}^{i-1}\geq x_j^i$.
\item $\gcn$ are the integer points of $\gc$.
\item Let $p:\gc\rightarrow \mathbb{R}^M$ be the projection that
picks out the last row of the triangular array, i.e. $p(x)=(x^M_1,\dots,x^M_M)$.
Likewise, let $q:\gcn\rightarrow \mathbb{N}^M$, the projection that
picks out the last row of an integer triangular array.
\end{enumerate}

\subsection{RSK}
Recall that a partition $\lambda$ of $k$ is a vector of
integers $(\lambda_1,\lambda_2,\dots)$, where $\lambda_1\geq \lambda_2
\geq\dots$ such that $\sum_i \lambda_i=k$. 
It follows that only finitely many of the $\lambda_i$:s are
non-zero.

A partition can be represented by a \emph{Young diagram},
drawn as a configuration of  boxes aligned in rows. 
The $i$:th row of boxes is $\lambda_i$ boxes long. 
A \emph{semi-standard Young tableau} (SSYT) is a filling of the 
boxes of a Young diagram with natural numbers, increasing from left to right
 in rows and strictly increasing from the top down  in columns. 
The \emph{Robinsson-Schensted-Knuth algorithm}  (RSK algorithm) is an algorithm
that bijectively maps $W_{M,N}$ to pairs of semi-standard Young tableau.
For details of this algorithm, see for example 
\cite{sagan:symgroup, stanley:ec2}.

Fix a matrix $w\in W_{M,N}$.
This matrix is mapped by RSK to a pair of SSYT,
$(P(w),Q(w))$. The $P$ tableau will contain elements of 
$\mathbf{M}:=\{1,2,\dots,M\}$ only. 
Construct a triangular array 
\begin{equation*}
x=\begin{matrix}
x_1^1\\
\vdots & \ddots\\
x_1^{M-1} & \dots & x_{M-1}^{M-1}\\
x_1^{M} & \dots & x_{M}^{M-1} &  x_{M}^{M}\\
\end{matrix}
\end{equation*}
where $x^i_j$ is the coordinate of the rightmost box 
filled with a number at most $i$ in the $j$-th row of the $P(w)$-tableau. 
This is a map from $W_{M,N}$ to $\gcn$.

\subsection{A Measure on Semistandard Young Tableau}
\label{sec:meas-semist-young}
Consider the following probability measure on $W_{M,N}$. The elements
in the matrix are i.i.d. geometric random variables with parameter $q^2$.
Recall that a  variable $X$ is geometrically distributed with parameter
$q^2$, written $X\in \Ge(q^2)$ if $P[X=k]=(1-q^2)(q^2)^{k}$,
$k\geq 0$. 
The square here will save a lot of root signs later. Such a
stochastic variable
has expectation $a=q^2/(1-q^2)$ and variance $b=q^2/(1-q^2)^2$.

Applying the RSK algorithm to this array
 induces a measure on SSYT:s, and by the correspondence above,
a measure on $\gcn$. Call this measure $\pi^{\RSK}_{q^2,M,N}$.
The following is shown in \cite{johansson:shapefluct}.
\begin{prop}
Let $W_{M,N}$ contain i.i.d. $\Ge(q^2)$ random variables in each position.
The probability that the RSK correspondence, when applied to this 
matrix, will yield Young diagrams of shape 
$\lambda=(\lambda_1,\dots,\lambda_M)$ is 
\begin{multline}
\rho^{\RSK}_{q^2,M,N}:=\frac{(1-q^2)^{MN}}{M!}
\prod_{j=0}^{M-1}
\frac{1}{j!(N-M+j)!} 
\times
\\
\times
\prod_{1\leq i<j\leq M}(\lambda_i-\lambda_j+j-i)^2 
\prod_{i=1}^M
\frac{(\lambda_i+1)!}{(\lambda_i+M-i)!} q^{2k},
\end{multline}
where $k=|\lambda|=\sum_i \lambda_i$.
\end{prop}

In other words, the measure $\pi^{\RSK}_{q^2,M,M}$, integrating out
all variables not on the last row, is $\rho^{\RSK}_{q^2,M,N}$.
This, together with the following result characterizes
the measure $\pi^{\RSK}_{q^2,M,M}$ completely.

\subsection{Uniform lift}
Proposition 3.2 in \cite{baryshnikov:gue} states that the 
probability measure $\pi^{\RSK}_{q^2,M,N}$,
conditioned on the last row of the triangular array being $\lambda$,
is uniform on the cone $q^{-1}(\lambda):=\{x\in \gcn: q(x)=\lambda\}$. 

In formulas, this can be formulated as follows.
\begin{prop}
\label{thm:rsklift}
For any bounded continuous function 
$\phi:\mathbf{M}\times\mathbb{Z}\rightarrow \mathbb{R}$ of compact support,
\begin{equation*}
\mathbb{E}_{\pi^{\RSK}_{q^2,M,N}}[\prod_{i,j} (1+\phi(i,x_j^i))]
= \sum_{\lambda} \left(
\frac{1}{L(\lambda)}
\sum_{x\in q^{-1}(\lambda)} 
\prod_{i,j} (1+\phi(i,x_j^i))\right)
\rho^{\RSK}_{q^2,M,N}(\lambda).
\end{equation*}
where $L(\lambda)$ is the number of integer points in 
$q^{-1}(\lambda)$.
\end{prop}
The number of such integer points is given
by 
\begin{equation*}
L(\lambda)=\prod_{i<j} \frac{\lambda_i-\lambda_j+j-i}{j-i}.
\end{equation*}
\subsection{GUE Eigenvalue measure}
It is well know, see for example \cite{mehta:randmat}, that 
the probability measure on the eigenvalues induced
by GUE measure on $M\times M$ hermitian matrices is
\begin{equation*}
\rho^{\GUE}_M(\lambda_1,\dots,\lambda_M)=
 \frac{1}{Z_M} 
\prod_{1\leq i<j\leq M} (\lambda_i-\lambda_j)^2
\prod_{1\leq i\leq M} \exp(-\lambda_i^2)
\end{equation*}
for some constant $Z_M$ that we need not be concerned with here.

\subsection{Uniform lift of GUE measure}
\cite{baryshnikov:gue}  shows a result for eigenvalues of 
minors of GUE matrices that is similar to the above result for partitions.
He shows that given the eigenvalues of the whole matrix
$\lambda=(\lambda_1>\dots> \lambda_M)$, the triangular array
of  eigenvalues of all the minors are uniformly distributed in 
$p^{-1}(\lambda):=\{x\in\gc:p(x)=\lambda\}$. 
Again we can write this more formally.
\begin{prop}
\label{thm:guelift}
For any bounded continuous function 
$\phi:\mathbf{M}\times\mathbb{R}\rightarrow\mathbb{R}$
 of compact support,
the measure $\pi^{\GUE}_M$ satisfies
\begin{equation*}
\mathbb{E}[\prod_{i,j}(1+\phi(i,\mu^i_j))]=
\int_\lambda
\left(
\frac{1}{\Vol(\lambda)}
\int_{p^{-1}(\lambda)} \prod_{i,j}(1+\phi(i,\mu^i_j))
\right) \rho^{\GUE}_M(\lambda)\, d^M\lambda.
\end{equation*}
where $\Vol(\lambda)$ is the volume of the cone~$p^{-1}(\lambda)$.
\end{prop}
This volume is given by
\begin{equation*}
\Vol(\lambda)=\prod_{i<j} \frac{\lambda_i-\lambda_j}{j-i}.
\end{equation*}
This situation is then very similar to the measure $\pi^{\RSK}_{q^2,M,N}$
above, in the sense that, conditioned the last row, the rest of the variables
is uniformly distributed in a certain cone.

\subsection{Scaling limit}
We are now in a position to see the connection between  
the measures $\pi^{\RSK}_{q^2,M,N}$ and $\pi^{\GUE}_M$.
\begin{prop}
\label{thm:rsktogue}
Let $a:=\mathbb{E}[w(1,1)]=q^2/(1-q^2)$ and 
$b:=\Var[w(1,1)]=q^2/(1-q^2)^2$.
Then for any bounded continuous function of compact support $\phi$,
\begin{equation*}
\mathbb{E}_{\pi^{\GUE}_M}[
\prod_{i,j}
(1+\phi(i,\mu_j^i))] = 
\lim_{N\rightarrow\infty}
\mathbb{E}_{\pi^{\RSK}_{q^2,M,N}}[
\prod_{i,j}
(1+\phi(i,\frac{x_j^i-aN}{\sqrt{b N}}))]. 
\end{equation*}
\end{prop}
\begin{proof}
Plug in the expression for the right hand side in 
 proposition~\ref{thm:rsklift} and for the left hand side in 
proposition~\ref{thm:guelift}.
Stirling's formula and the convergence of a Riemann sum to an
integral proves the theorem.
\end{proof}

\subsection{Polynuclear growth}
\label{sec:polynuclear-growth}
The measure $\pi^{\RSK}_{q^2,N,M}$ 
is a version of the Schur process
and is a determinantal process on $\mathbf{M}\times \mathbb{N}$,
by~\cite{okounkov:schur}.
We will use the following result from 
\cite{johansson:detproc}. 
\begin{prop}
\label{thm:png}
The process $\{x_j^i\}$ with the measure described 
in~\ref{sec:meas-semist-young} is determinantal with kernel
\begin{equation}
\label{eq:4}
K^{\PNG}_{q^2,N,M}(r, x, s, y)=
\frac{1}{(2\pi i)^2}
\int\frac{dz}{z} \int\frac{dw}{w} \frac{z}{z-w}
\frac{w^{y+N}}{z^{x+N}} 
\frac{(1-qw)^s}{(1-qz)^r}
\frac{(z-q)^{N-M}}{(w-q)^{N-M}}.
\end{equation}
For $r\leq s$, the paths of integration 
for  $z$ and $w$ are anticlockwise 
along 
circles centred at zero with radii such that $q< |w| <  |z| < 1/q$.
For the case $r > s$, integrate instead along circles 
such that $q<|z|<|w|<1/q$. 
\end{prop}
This follows immediately from proposition 3.12 and theorem 3.14 in 
\cite{johansson:detproc}. 

Having now introduced the $\PNG$-kernel,
 we can state the following scaling limit result.
\begin{lemma}
\label{thm:pnglimit}
Let $a=q^2/(1-q^2)$ and $b=q^2/(1-q^2)^2$ as above.
The following claims are true for $M$ fixed.
\begin{enumerate}
\item
For $r$, $s\leq M$, 
 \begin{equation*}
\frac{g(r, \xi,N)}{g(s,\eta,N)}\sqrt{2bN}
K^{\PNG}_{N,M}(r, \lfloor aN+\xi\sqrt{2bN}\rfloor; \,
s, \lfloor aN+\eta\sqrt{2bN}\rfloor)\longrightarrow
K^{\Herm}(r, \xi;\, s, \eta)
\end{equation*}
 uniformly on compact sets
as $N\rightarrow\infty$ 
 for a certain function $g\neq 0$.
\item
The expression
\begin{equation*}
\frac{g(r, \xi,N)}{g(s,\eta,N)}
\sqrt{2bN}K^{\PNG}_{N,M}(r, \lfloor aN+\xi\sqrt{2bN}\rfloor;\, s, \lfloor aN+\eta\sqrt{2bN}\rfloor)
\end{equation*}
is bounded uniformly for $1\leq r,s\leq M$
and $\xi$, $\eta$ in a compact set.
\end{enumerate}
\end{lemma}

The proof, given in section \ref{sec:calc}, 
is an asymptotic analysis of the integral in~(\ref{eq:4}).
Now everything is set up so we can prove the main result of this section.
\begin{proof}[Proof of theorem \ref{thm:gue}.]
According to proposition~\ref{thm:rsktogue},  
\begin{equation}
\label{eq:16}
\mathbb{E}_{\pi^{\GUE}_M}[\prod(1+\phi(i,\mu^i_j))]=
\lim_{N\rightarrow\infty}
\mathbb{E}_{\pi^{\RSK}_{q^2,M,N}}[
\prod_{i,j}
(1+\phi(i,\frac{x_j^i-aN}{\sqrt{b N}}))] .
\end{equation}
The point processes on the right hand side of this last
expression are determinantal. 
Their kernels can be written 
\begin{equation*}
K^N(r,\xi,s,\eta):=\frac{g(r, \xi,N)}{g(s,\eta,N)}\sqrt{2bN}
K^{\PNG}_{N,M}(r, aN+\xi\sqrt{2bN}; \,
s, aN+\eta\sqrt{2bN}).
\end{equation*}
for some function $g$ that cancels out in all determinants, and therfore
does not affect the correlation functions.
By lemma~\ref{thm:pnglimit}, these $K^N$ satisfy all the assumptions of 
proposition~\ref{thm:determ-point-proc}.
Thus, the point processes that these define converge weakly
to a point process with kernel $K^{\Herm}$.
This implies that the measure on the left hand side of 
equation~(\ref{eq:16}), i.e. $\pi^{\GUE}_M$ is determinantal 
with kernel $K^{\Herm}$. The observation that  $M$ was arbitrary
completes the proof. 
\end{proof}

\section{Aztec Diamond}
\label{sec:aztec-diamond-1}The point-process connected to the tilings of 
this shape, described in the introduction was thoroughly
analyzed in \cite{johansson:aztec}. 
The following result is shown.
\begin{prop}
\label{thm:aztec}
The process $\{x_j^i\}$ described in section~\ref{sec:aztec-diamond} 
is determinantal on
$\Lambda=\mathbb{N}\times\mathbb{N}$, with kernel 
$K^A_N$ given by
\begin{equation}
K^{\Aztec}_N(2r,x,2s,y)=
\frac{1}{(2\pi i)^2}
\int\frac{dz}{z} \int \frac{dw}{w}
\frac{w^y(1-w)^s (1+1/w)^{N-s}}
{z^x(1-z)^r(1+1/z)^{N-r}}
\frac{z}{z-w}
\end{equation}
and reference measure $\mu$ which is counting measure on $\mathbb{N}$.
The paths of integration are as follows: For $r\leq s$,
integrate $w$ along a contour enclosing its pole at $-1$ anticlockwise,
and $z$ along a contour enclosing $w$ and the pole at $0$ but not the 
pole at $1$ anticlockwise. For $r>s$, switch the contours of 
$z$ and $w$.  
\end{prop}
Based on this integral formula we can prove
the following scaling limit analogous to that in lemma~\ref{thm:pnglimit}.

\begin{lemma}
\label{krawtchouk_to_hermite}
The following claims hold.
\begin{enumerate}
\item[(1)]
\begin{equation*}
\frac{g(r,\xi,N)}{g(s,\eta,N)}
\sqrt{N/2}
K_N^{\Aztec}(2r, \lfloor N/2+\xi\sqrt{N/2}\rfloor;\,2s,\lfloor N/2+\eta\sqrt{N/2}\rfloor)
\longrightarrow K^{\Herm}(r,\xi;\,s,\eta)
\end{equation*}
uniformly on compact sets as $N\rightarrow\infty$ 
for an appropriate  function $g\neq 0$.
\item[(2)] The expression
\begin{equation*}
\frac{g(r,\xi,N)}{g(s,\eta,N)}
\sqrt{N/2}
K_N^{\Aztec}(2r, \lfloor N/2+\xi\sqrt{N/2}\rfloor;
\,2s,\lfloor N/2+\eta\sqrt{N/2}\rfloor)
\end{equation*}
is uniformly bounded with respect to $N$ for  $(r$, $\xi$, $s$, $\eta)$ contained
in any  compact set in $\mathbb{N}\times \mathbb{R}\times \mathbb{N}\times \mathbb{R}$.
\end{enumerate} 
\end{lemma}
The proof is based on a saddle point
  analysis that is presented it section~\ref{sec:calc}.
We can now set about proving the main result of this section.

\begin{proof}[Proof of theorem \ref{thm:main_aztec_thm}]
By proposition~\ref{thm:aztec}, the $x^i_j$ form a determinantal process
with kernel $K^{\Aztec}_N$. The rescaled process $(x^i_j-N/2)/\sqrt{N/2}$
has kernel
\begin{equation}
K^N(r,\xi; s,\eta):=\frac{g(r,\xi,N)}{g(s,\eta,N)}\sqrt{N/2}K_N^{\Aztec}(2r,N/2+\xi\sqrt{N/2};\,2s, N/2+\eta\sqrt{N/2}).
\end{equation}
By lemma~\ref{krawtchouk_to_hermite}, the kernels $K^N$ satisfy
all the assumptions of proposition~\ref{thm:determ-point-proc}.
So they converge to the process with kernel $K^{\GUE}$.
\end{proof}

\section{The Hexagon}
Consider an (a,b,c)-hexagon, such as the one 
in figure~\ref{fig:small_hexagon}. 
%We will show that the blue dots form a determinantal process. 
First we need some coordinate system to describe the position of the dots.
Say that the $a$ simple, symmetric, random walks start at $t=0$ and
$y=0$, $2$, \dots, $2a-2$. 
In each unit of time, they move one unit up or down,
and are conditioned to end up at $y=c-b$, $c-b+2$, \dots, $c-b+2a-2$ at
time $t=b+c$ and never to intersect.
One realisation of this process is the red dots in 
figure~\ref{fig:small_hexagon}. 
At time $t$,  the only possible $y$-coordinates for the 
red dots are $\{\alpha_t+2k\}_{0\leq k \leq \gamma_t}$, 
where
\begin{align*}
\gamma_t&=\begin{cases}
t+a-1 & 0\leq t\leq b\\
b+a-1 & b\leq t\leq c\\
a+b+c-t-1 &c\leq t\leq b+c,
\end{cases} &
\alpha_t&=\begin{cases}
-t & 0\leq t\leq b\\
t-2b &b\leq t\leq b+c.
\end{cases}
\end{align*}
Let $\Lambda_{a,b,c}=\{ (t, \alpha_t+2k):0\leq t\leq b+c,\, 0\leq k\leq \gamma_t\} $ 
be the set of all the dots, red and blue. 

\subsection{A determinantal kernel for the hexagon tiling process}
We now need to define the normalised associated Hahn polynomials,
$\tilde q_{n,N}^{(\alpha,\beta)}(x)$. These
orthogonal polynomials satisfy
\begin{equation}
\label{eq:hahn_orthogonalityrelation}
\sum_{x=0}^N 
\tilde q_{n,N}^{(\alpha,\beta)}(x)
\tilde q_{m,N}^{(\alpha,\beta)}(x)
\tilde w_N^{(\alpha,\beta)}(x)  =\delta_{n,m},
\end{equation}
where the weight function is 
\begin{equation*}
\tilde w_N^{(\alpha,\beta)}(x) 
=\frac{1}{x!(x+\alpha)!(N+\beta-x)!(N-x)!}.
\end{equation*}
They can be computed as
\begin{equation*}
\tilde q^{(\alpha,\beta)}_{n,N}(x)=
\frac{(-N-\beta)_n(-N)_n}{\tilde d_{n,N}^{(\alpha,\beta)}n!}
\,_3F_2(\substack{-n,n-2N-\alpha-\beta-1, -x\\-N-\beta,-N};1),
\end{equation*}
where
\begin{equation*}
\left(\tilde d_{n,N}^{(\alpha,\beta)}\right)^2=
\frac{(\alpha+\beta+N-1-n)_{N+1}}
{(\alpha+\beta+2N+1-2n)n!(\beta+N-n)!(\alpha+N-n)!(N-n)!}.
\end{equation*}

For convenience, let $a_r=|c-r|$
and $b_r=|b-r|$.
\cite{johansson:rhombus} shows the following.
\begin{prop}
\label{thm:red-dots}
The red dots form a determinantal point process
on the space $\Lambda_{a,b,c}$ with kernel
\begin{multline*}
\tilde K^{\Loz}_{a,b,c} (r, \alpha_r + 2x;\,
s,\alpha_s+2y)=
-\phi_{r,s}(\alpha_r + 2x,\alpha_s+2y)\\
+\sum_{n=0}^{a-1}
\sqrt{\frac{(a+s-1-n)!(a+b+c-r-1-n)!}{(a+s-1-n)! (a+b+c-s-1-n)!}}
\tilde q_{n,\gamma_r}^{b_r,a_r} (x)
\tilde q_{n,\gamma_s}^{b_s,a_s} (y)
\omega_r(x) \tilde \omega_s(y),  
\end{multline*}
where $\phi_{r,s}(x,y)=0$ if $r\geq s$ and 
\begin{equation*}
\phi_{r,s}(x,y)=\binom{s-r}{\frac{y-x+s-r}{2}}
\end{equation*}
otherwise. Furthermore, 
\begin{align*}
\omega_r(x)&=\begin{cases}
((b_r+x)!(\gamma_r+a_r-x)!)^{-1} & 0\leq r\leq b\\
(x!(\gamma_r+a_r-x)!)^{-1} & b\leq r\leq c\\
(x!(\gamma_r-x)!)^{-1} & c\leq r\leq b+c
\end{cases}\\
\intertext{and}
\tilde \omega_s(y)&=\begin{cases}
(y!(\gamma_s-y)!)^{-1} 
& 0\leq r\leq b\\
((b_s+y)!(\gamma_s-y)!)^{-1}
& b \leq r \leq c\\
((b_r+y)!(\gamma_r+a_r-y)!)^{-1} 
& c\leq r\leq b+c.
\end{cases}
\end{align*}
\end{prop}

It follows that the blue dots also form a determinantal point process. 
To compute its kernel we need the following lemma.
\begin{lemma}
\label{thm:hexagon_sum}
\begin{multline*}
\binom{s-r}{\frac{s-r+2y+\alpha_s-2x-\alpha_r}{2}}=\\
\sum_{n=0}^\infty 
\sqrt{\frac{(a+s-1-n)!(a+b+c-r-1-n)!}{(a+r-1-n)!(a+b+c-s-1-n)!}}
\tilde q_{n,\gamma_r}^{(b_r,a_r)}(x)
\tilde q_{n,\gamma_s}^{(b_s,a_s)}(y)
\omega_r(x)
\tilde\omega_s(y)
\end{multline*}
when $ s\geq r$.
\end{lemma}
\begin{proof}
This proof uses the results obtained in the proof of~\ref{thm:red-dots} 
in~\cite[equation 3.25]{johansson:rhombus}. 
Define convolution product as follows. For 
$f, g:\mathbb{Z}^2\rightarrow\mathbb{Z}$, define $(f*g):\mathbb{Z}^2\rightarrow\mathbb{Z}$ by
\begin{equation*}
(f*g)(x,y):= \sum_{z\in\mathbb{Z}} f(x,z) g(z,y).
\end{equation*}
Let $\phi(x,y):=\delta_{x,y+1}+ \delta_{x,y-1}$. 
Also let
\begin{align*}
\phi^{*0}(x,y)&:=\delta_{x,y}\\
\phi^{*1}(x,y)&:=\phi(x,y)\\
\phi^{*n}(x,y)&:=(\phi^{*(n-1)}*\phi)(x,y).
\end{align*}
Set
\begin{align*}
c_{j,k}&:=\frac{1}{(a-k)(j-k)!(a-1-j)!}\\
f_{n,k}&:=\binom{n}{k}\frac{(n-2a-b-c+1)_k}{(-a-b+1)_k(-a)_k}
\end{align*}
and finally let
\begin{align*}
\psi(n,z)&:=\sum_{m=0}^n f_{n,m} \sum_{j=m}^{a-1} c_{j,m}\phi(2j,z)\\
\phi_{0,1}(n,y)&:=\psi(n,y)\\
\phi_{0,r}(n,y)&:=\psi*\phi^{*(r-1)}(n,y).
\end{align*}
The dual orthogonality relation to~(\ref{eq:hahn_orthogonalityrelation}) 
is precisely
\begin{equation}
\label{eq:8}
\sum_{n=0}^{\gamma_r}
\tilde q_{n,\gamma_r}^{(b_r,a_r)}(x)
\tilde q_{n,\gamma_r}^{(b_r,a_r)}(y)
\omega_r(x)\tilde \omega_r(y)=\delta_{x,y}.
\end{equation}

By equation (3.25), (3.30) and (3.32) of the above mentioned paper,
\begin{equation*}
\label{eqn:Aabcrn}
\phi_{0,r}(n,\alpha_r+2z)=
 A(a,b,c,r,n)
\tilde q_{n,\gamma_r}^{(b_r,a_r)}(z)
\tilde \omega_r(z),
\end{equation*}
where 
\begin{equation*}
A(a,b,c,r,n):=\frac{(a+1)_{r-1} \tilde d_{n,\gamma_r}^{(b_r,a_r)}n!}{
(-a-c+1)_n(-a-b-c+r+1)_n}.
\end{equation*}
Inserting this into the orthogonality relation in~(\ref{eq:8})  gives
\begin{equation*}
\sum_{n=0}^{\gamma_r}
\tilde q_{n,\gamma_r}^{(b_r,a_r)}(x)
\frac{\omega_r(x)}{A(a,b,c,r,n)}
\phi_{0,r}(n,\alpha_r+2z)=\delta_{x,y}.
\end{equation*}
Convolving both sides of the above relation with $\phi^{*(s-r)}$ gives
\begin{multline*}
\sum_{n=0}^{\gamma_r}
\tilde q_{n,\gamma_r}^{(b_r,a_r)}(x)
\frac{\omega_r(x)}{A(a,b,c,r,n)}
\sum_{z\in \mathbb{Z}}
\phi_{0,r}(n,\alpha_r+2z)
\phi^{*(s-r)}(\alpha_r+2z, \alpha_s+2y)
\\=\phi^{*(s-r)}(\alpha_r+2x, \alpha_s+2y),
\end{multline*}
which, when the left hand side is simplified, gives
\begin{equation*}
\sum_{n=0}^{\gamma_r}
\tilde q_{n,\gamma_r}^{(b_r,a_r)}(x)
\frac{\omega_r(x)}{A(a,b,c,r,n)}
\phi_{0,s}(n,\alpha_r+2y)
=\phi^{*(s-r)}(\alpha_r+2x, \alpha_s+2y).
\end{equation*}
Invoking equation~(\ref{eqn:Aabcrn}) again to simplify the left hand side
 and explicitly calculating
the right hand side gives 
\begin{equation*}
\sum_{n=0}^{\gamma_r}
\frac{A(a,b,c,s,n)}{A(a,b,c,r,n)}
\tilde q_{n,\gamma_r}^{(b_r,a_r)}(x)
\tilde q_{n,\gamma_s}^{(b_s,a_s)}(y)
\omega_r(x)\tilde\omega_s(y)=
\binom{s-r}{\frac{s-r+2y+\alpha_s-2x-\alpha_r}{2}}.
\end{equation*}
It is easy to check that 
\begin{equation*}
\frac{A(a,b,c,s,n)}{A(a,b,c,r,n)}=
\sqrt{\frac{(a+s-1-n)!(a+b+c-r-1-n)!}{(a+r-1-n)!(a+b+c-s-1-n)!}},
\end{equation*}
which proves the lemma.
\end{proof}

We now need to introduce the normalized Hahn polynomials 
$q_{n,N}^{(\alpha,\beta)}(x)$.
These satisfy
\begin{equation}
\sum_{x=0}^N
q^{(\alpha,\beta)}_{n,N}(x) q^{(\alpha,\beta)}_{m,N}(x) w^{(\alpha,\beta)}_N(x)
=\delta_{m,n},
\end{equation}
where
\begin{equation}
w^{(\alpha,\beta)}_N(t) = \frac{(N+\alpha-t)!(\beta+t)!}{t!(N-t)!}.
\end{equation}
\begin{thm}
\label{thm:blue-dots}
The blue dots form a determinantal point process on the space 
$\Lambda_{a,b,c}$ with kernel 
\begin{multline*}
K^{\Loz}_{a,b,c}(r,x;\,s,y) = \\
\sum_{n=-\infty}^{-1}
\sqrt{\frac{(s+n)!(b+c-r+n)!}{(r+n)!(b+c-s+n)!}}
 q_{r+n,\gamma_r}^{(b_r,a_r)}(x)
 q_{s+n,\gamma_s}^{(b_s,a_s)}(y)
\sqrt{w^{(b_r,a_r)}_{\gamma_r}(x)w^{(b_s,a_s)}_{\gamma_s}(y)},
\end{multline*}
when $s\geq r$, and
\begin{multline*}
K^{\Loz}_{a,b,c}(r,x;\,s,y) =\\ -
\sum_{n=0}^{a-1}
\sqrt{\frac{(s+n)!(b+c-r+n)!}{(r+n)!(b+c-s+n)!}}
 q_{r+n,\gamma_r}^{(b_r,a_r)}(x)
 q_{s+n,\gamma_s}^{(b_s,a_s)}(y)
\sqrt{w^{(b_r,a_r)}_{\gamma_r}(x)w^{(b_s,a_s)}_{\gamma_s}(y)}
\end{multline*}
otherwise.
\end{thm}
\begin{proof}
It is well known that the complement of a  determinantal point processes on 
a finite set with kernel $K$  is also
determinantal with kernel $\tilde K=I-K$, i.e. 
$\tilde K(x,y)=\delta_{x,y}-K(x,y)$.

Applying this result to our problem, we consider 
$\delta_{x,y}\delta_{r,s}-\tilde K^{\Loz}_{a,b,c}(r,x;\,s,y)$.
We now separate two cases. When $s\geq r$ see  that
\begin{multline*}
K(r,x;\,s,y) =  \binom{s-r}{y-x+\frac{s-r+\alpha_s-\alpha_r}{2}}-\\
\sum_{n=0}^{a-1}
\sqrt{\frac{(a+s-1-n)!(a+b+c-r-1-n)!}{(a+r-1-n)!(a+b+c-s-1-n)!}}
\tilde q_{n,\gamma_r}^{(b_r,a_r)}(x)
\tilde q_{n,\gamma_s}^{(b_s,a_s)}(y)
\omega_r(x)
\tilde\omega_s(y)
\end{multline*}
is a candidate for the kernel for the blue particles.
By lemma \ref{thm:hexagon_sum} this simplifies to 
\begin{multline*}
K(r,x;\,s,y) = \\
\sum_{n=a}^{\infty}
\sqrt{\frac{(a+s-1-n)!(a+b+c-r-1-n)!}{(a+r-1-n)!(a+b+c-s-1-n)!}}
\tilde q_{n,\gamma_r}^{(b_r,a_r)}(x)
\tilde q_{n,\gamma_s}^{(b_s,a_s)}(y)
\omega_r(x)
\tilde\omega_s(y).
\end{multline*}
For $s<r$ we just get
\begin{multline*}
K(r,x;\,s,y) = \\-
\sum_{n=0}^{a-1}
\sqrt{\frac{(a+s-1-n)!(a+b+c-r-1-n)!}{(a+r-1-n)!(a+b+c-s-1-n)!}}
\tilde q_{n,\gamma_r}^{(b_r,a_r)}(x)
\tilde q_{n,\gamma_s}^{(b_s,a_s)}(y)
\omega_r(x)
\tilde\omega_s(y).
\end{multline*}

We now exploit a useful duality result from
\cite{borodin:orthogonality_on_a_finite_set}. 
It  states that
\begin{equation*}
q^{(\alpha,\beta)}_{n,N}(x)\sqrt{w^{(\alpha,\beta)}_N(x)} = 
(-1)^{x}\tilde q^{(\alpha,\beta)}_{N-n,N}(x)
\sqrt{\tilde w^{(\alpha,\beta)}_N(x)}.
\end{equation*}

Insert this into the formulas above and define the new kernel
\[
K^{\Loz}_{a,b,c}(r,x;\,s,y):= 
(-1)^{y-x}\sqrt{\omega_s(y)\omega_r(x)^{-1}
\tilde \omega_r(x)\tilde \omega_s(y)^{-1}} K(r,x;\,s,y).
\]
This kernel gives the same correlation functions as $K$, since the extra
factors cancel out in the determinants. The new kernel can be written as
\begin{multline*}
K^{\Loz}_{a,b,c}(r,x;\,s,y) = \sqrt{w^{(b_r,a_r)}_{\gamma_r}(x)w^{(b_s,a_s)}_{\gamma_s}(y)} \times \\
\sum_{n=a}^{\infty}
\sqrt{\frac{(a+s-1-n)!(a+b+c-r-1-n)!}{(a+r-1-n)!(a+b+c-s-1-n)!}}
 q_{\gamma_r-n,\gamma_r}^{(b_r,a_r)}(x)
 q_{\gamma_s-n,\gamma_s}^{(b_s,a_s)}(y),
\end{multline*}
when $s\geq r$, and
\begin{multline*}
K^{\Loz}_{a,b,c}(r,x;\,s,y) =-\sqrt{w^{(b_r,a_r)}_{\gamma_r}(x)w^{(b_s,a_s)}_{\gamma_s}(y)} \times \\ 
\sum_{n=0}^{a-1}
\sqrt{\frac{(a+s-1-n)!(a+b+c-r-1-n)!}{(a+r-1-n)!(a+b+c-s-1-n)!}}
 q_{\gamma_r-n,\gamma_r}^{(b_r,a_r)}(x)
 q_{\gamma_s-n,\gamma_s}^{(b_s,a_s)}(y)
\end{multline*}
otherwise.

The change of variables $j:=a-1-n$ puts these expressions 
on a simpler form, thereby proving the theorem.
\end{proof}

\subsection{Asymptotics}
Let $0<p<1$ be some real number. Let $\alpha=\gamma pN$, $\beta=\gamma(1-p)N$,
$\tilde x= \lfloor pN+\sqrt{2p(1-p)N(1+\gamma^{-1})} x\rfloor$. Then
\begin{equation}
\label{eq:13}
\sqrt[4]{2p(1-p)N(1+\gamma^{-1})}
\sqrt{w_{n,N}^{(\alpha,\beta)}(\tilde x)} q_{n,N}^{(\alpha,\beta)}(\tilde x)
\longrightarrow 
(-1)^n \sqrt{e^{-x^2}} h_n(x)
\end{equation}
uniformly on compact sets in $x$ as $N\rightarrow\infty$.
For completeness we give the proof of this result in the appendix.

We would like to apply this with  $p=\frac{1}{2}$ 
and $\gamma=2$ to our kernel $K^{\Loz}$ 
and letting $a=b=c\rightarrow\infty$, i.e. we would like to take
the limit
\begin{multline*}
K(r,\xi; s,\eta) =\\
 \lim_{N=a=b=c\rightarrow\infty} (-3N)^{r-s}\sqrt{3N/4}K^{\Loz}_{a,b,c} 
(r, \lfloor N/2+\xi \sqrt{3N/4}\rfloor;
s, \lfloor N/2+\eta \sqrt{3N/4}\rfloor).
\end{multline*}
The factor $(-3N)^{r-s}$ cancels out in all determinants and 
is thus of no import.
For $s\geq r$ we get
\begin{equation}
K(r,\xi; s,\eta) = 
\sum_{j=-\infty}^{-1}\sqrt{\frac{(s+j)!}{(r+j)!}}
h_{r+j}(\xi)h_{s+j}(\eta) e^{-(\xi^2+\eta^2)/2}
\end{equation}
and formally, if we ignore the fact that this turns into 
an infinite sum, for $s<r$ we get
\begin{equation}
\label{eq:14}
K(r,\xi; s,\eta) = 
-\sum_{j=0}^{\infty}\sqrt{\frac{(s+j)!}{(r+j)!}}
h_{r+j}(\xi)h_{s+j}(\eta) e^{-(\xi^2+\eta^2)/2}.
\end{equation}
This expression can be simplified with the following lemma
\begin{lemma}
\label{thm:hermitesum}
Let $H$ be the Heaviside function defined by equation~(\ref{eq:heaviside})
above.
Then,
\begin{multline}
\label{eq:2}
\frac{\sqrt{2^k}}{(k-1)!}(x-y)^{k-1} H(x-y)= 
\sum_{n=k}^\infty
\sqrt{\frac{(n-k)!}{n!}}
h_{n-k}(y)h_{n}(x) e^{-y^2}\\
+\frac{1}{\sqrt[4]{\pi}}\sum_{n=0}^{k-1} 
\frac{h_n(x)\sqrt{2^{k-n}}}{\sqrt{n!}(k-1-n)!}
\int_{y}^\infty (t-y)^{k-1-n} e^{-t^2} dt
\end{multline}
pointwise for $x\neq y$.
\end{lemma}
The proof is given in section \ref{sec:calc}.

In view of this result, the  infinite series in  \ref{eq:14} converges
and the kernel $K$ is exactly the GUE minor kernel $K^{\Herm}$. 
The interpretation of this is the following. 
The distribution of the blue particles,
properly rescaled, tends weakly to the distribution of the eigenvalues
of GUE minors as the size of the diamond tends to infinity, 
equation~(\ref{eq:15}).
The only thing needed to make this a theorem is appropriate estimates
of the Hahn polynomials to control the 
convergence to the infinite sum.

\section{Proof of lemmas}
\label{sec:calc}

\begin{proof}[Proof of lemma \ref{thm:hermitesum}]
As the Hermite polynomials are orthogonal,
there is an expansion of the function in the left hand side
of~(\ref{eq:2}) of the form
\begin{equation}
\label{eq:3}
(x-y)^{k-1}H(x-y)=\sum_{n=0}^\infty
c_n(y) H_n(x),
\end{equation}
where $H_n$ is the $n$:th Hermite polynomial, as defined in
for example \cite{koekoek:askey}, and
where the coefficients are given by
\begin{equation}
c_n(y)=\frac{1}{2^{n}n!\sqrt{\pi}} 
\int_{-\infty}^\infty (x-y)^{k-1}H(x-y) H_n(x)e^{-x^2}dx.
\end{equation}
It is known that $e^{-x^2} H_n(x)=-\frac{d}{dx}(e^{-x^2}H_{n-1}(x))$
for $n\geq1$. 
Integration by parts and limiting the integration interval
according to the Heaviside function gives

\begin{equation*}
\int_{y}^\infty (x-y)^{k-1} H_n(x)dx=
\int_y^\infty (k-1)(x-y)^{k-2} H_{n-1}(x)dx
\end{equation*}

For $n\geq k$, repeat this process $k-1$ times to get
\begin{equation}
\label{eq:5}
c_n(y)=\frac{ (k-1)! e^{-y^2}H_{n-k}(y)}{2^{n}n!\sqrt{\pi}}.
\end{equation}
For $0\leq n <k$, stop doing partial integrations when $H_0$ 
is reached, giving
\begin{equation}
\label{eq:6}
\frac{(k-1)!}{(k-n-1)!}\int_y^\infty (x-y)^{k-n-1} e^{-x^2} \, dx.
\end{equation}
Inserting~(\ref{eq:5}) and~(\ref{eq:6}) in (\ref{eq:3}) and changing
to normalized Hermite polynomials proves the lemma. 
\end{proof}

\newcommand{\rvariable}{r}
\newcommand{\svariable}{s}
\begin{proof}[Proof of lemma~\ref{thm:pnglimit}.]
Assume first that $r\leq s$. By proposition~\ref{thm:png}
we have to consider the integral
\begin{equation*}
\frac{1}{(2\pi i)^2}
\int_{\gamma_{r_2}} dz
\int_{\gamma_{r_1}} \frac{dw}{w}
\frac{1}{z-w}
e^{Nf(z)-Nf(w)}
\frac{w^{\eta\sqrt{2bN}}}{z^{\eta\sqrt{2bN}}}
\frac{(1-qw)^2}{(1-qz)^r}
\frac{(w-q)^M}{(z-q)^M},
\end{equation*}
where $\gamma_r$ is a circle around the origin 
with radius $r$ oriented anticlockwise, $q<r_1<r_2<1/q$,
and 
\begin{equation}
f(z)=\log(z-q)-(1-q^2)^{-1} \log z.
\end{equation}
(Here we have ignored the difference between $aN+\xi\sqrt{2bN}$
and its integer part.)
Note that $f'(z)=0$ gives $z=1/q$. This leads us to choose
\begin{equation*}
r_1=\frac{1}{q}-\frac{2}{a\sqrt{N/2}},
\end{equation*}
and to deform $\gamma_{r_2}$ to a circle $\Gamma$
oriented clockwise
around
$1/q$ with radius $1/a\sqrt{N/2}$. 
The specific choice of radii are convenient for the 
computations below.
Choose
\begin{equation*}
g(r,\xi,N)=
2^{-r/2}
e^{-\xi^2/2}
q^{-\xi\sqrt{2bN}}
\left(\frac{q}{a\sqrt{N/2}}\right)^r.
\end{equation*}
Then, 
\begin{multline}
\label{eq:k1}
\frac{g(r,\xi,N)}{g(s,\eta,N)}
\sqrt{2bN}
K^{\PNG}_{N,M}(r, \lfloor aN+\xi\sqrt{2bN}\rfloor;\, s, \lfloor aN+\eta\sqrt{2bN}\rfloor)
\\
=\sqrt{2^{s-r}e^{\eta^2-\xi^2} }
q^{(\eta-\xi)\sqrt{2bN}}
\left(\frac{q}{a\sqrt{N/2}}\right)^{r-s}
\frac{\sqrt{2bN}}{(2\pi i)^2}
\\
\times
\int_{\Gamma} dz \int_{\gamma_{r_1}}\frac{dw}{w}
\frac{1}{z-w}
e^{Nf(z)-Nf(w)}
\frac{w^{\eta\sqrt{2bN}}}{z^{\eta\sqrt{2bN}}}
\frac{(1-qw)^2}{(1-qz)^r}
\frac{(w-q)^M}{(z-q)^M}.
\end{multline}
Parameterize $\gamma_{r_1}$ by $w(t)=r_1 e^{itE_n}$, 
$-\pi/E_N\leq t\leq \pi/E_N$,
$E_N=q/a\sqrt{N/2}$. We have 
\begin{align*}
\Real (f(w(0))-f(w(t))) & = \ln\left |\frac{w(0)-q}{w(t)-q}\right|\\
&=-\frac{1}{2} \ln\left(1+\frac{2r_1q(1-\cos E_Nt)}{(r_1-q)^2}\right)\\
&\leq -\frac{1}{2} \ln\left(1+q^2(1-\cos E_N t)\right),
\intertext{for $N$ large enough. Since $\cos x\leq 1+x^2/8$ when
$|x|\leq \pi$, the last expression is }
&\leq-\frac{1}{2} \ln\left(1+q^2 E_N^2 t^2/8)\right)\leq -Ct^2/N
\end{align*}
for $|t|\leq \pi/E_N$, where $C>0$ is a constant depending only on $q$.
Hence,
\begin{equation}
\label{eq:k2}
\Real N(f(w(0))-f(w(t))) \leq -Ct^2
\end{equation}
for $|t|\leq \pi/E_n$ with  $C>0$.

In  the right hand side of~(\ref{eq:k1}) we make the change of variables
\begin{equation}
\label{eq:k3}
z=z(u)=1/q-u/a\sqrt{N/2}
\end{equation}
with $u$ on the unit circle oriented anticlockwise. 
We obtain the integral 
\begin{multline}
\label{eq:k4}
\sqrt{2^{s-r}e^{\eta^2-\xi^2}}
\frac{2iq\sqrt{b}}{(2\pi i)^2 a}
\int_{\gamma_1} du \int_{-\pi/E_N}^{\pi/E_N} dt
\frac{1}{a\sqrt{N/2}(z(u)-w(t))}
e^{N(f(z(u))-f(w(t)))}\\
\times
\frac{(qw(t))^{\eta\sqrt{2bN}}}{(qz(u))^{\xi\sqrt{2bN}}}
\left(\frac{q}{a\sqrt{N/2}}\right)^{r-s}
\frac{(1-qw(t))^s}{(1-qz(u))^r}\frac{(w(t)-q)^M}{(z(u)-q)^M}.
\end{multline}
Note that $q\sqrt{b}/a=1$. Also, 
\begin{equation}
\label{eq:k4p}
f(1/q+h)=f(1/q)-a^2h^2/2+O(h^3)
\end{equation}
for $|h|$ small. Hence, for $N$ sufficiently large,
\begin{equation}
\label{eq:k5}
N(f(z(u))-f(w(0)))=-u^2+4+h_N(u)/\sqrt{N},
\end{equation}
where $h_N(u)$ is bounded for $|u|=1$.
We have
\begin{equation}
\label{eq:k6}
\frac{(qw(t))^{\eta\sqrt{2bN}}}{(qz(t))^{\xi\sqrt{2bN}}}
=\left(1-\frac{2q}{a\sqrt{N/2}} \right)^{\eta\sqrt{2bN}}
e^{2i\eta t}
\left(1-\frac{qu}{a\sqrt{a\sqrt{N/2}}}\right)^{-\xi\sqrt{2bN}}.
\end{equation}
By the inequality $(1+x/n)^n\leq e^x$ for $x>-n$, $n\geq 1$, 
the right hand side in~(\ref{eq:k5}) has a bound independent of $N$. 
We also have 
\begin{equation}
\label{eq:k7}
a\sqrt{N/2} |z(u)-w(t)|\geq 1
\end{equation}
for $u\in \gamma_1$, $|t|\leq \pi/E_N$, and 
\begin{equation}
\label{eq:k8}
\left|\left(\frac{q}{a\sqrt{N/2}}\right)
\frac{(1-qw(t))^s}{(1-qz(t))^r}\frac{(w(t)-q)^M}{(z(u)-q)^M}\right|
\leq CN^{s/2}
\end{equation}
for $u\in \gamma_1$, $|t|\leq \pi/E_N$, by~(\ref{eq:k3})
and the definition of $w(t)$.

It follows from~(\ref{eq:k2}),~(\ref{eq:k5}),~(\ref{eq:k6}),~(\ref{eq:k7})
and~(\ref{eq:k8}) that the part of the integral in~(\ref{eq:k4})
where the t-integration is restricted to $N^{1/3}\leq |t|\leq \pi/E_N$
can be bounded by
\begin{equation*}
CN^{s/2}
\int_{|t|\geq N^{1/3}}
e^{-Ct^2}\,dt,
\end{equation*}
which goes to $0$ as $N\rightarrow\infty$. When $|t|\leq N^{1/3}$
we have 
\begin{equation}
\label{eq:k9}
\left| \frac{a\sqrt{N/2}}{q} (1-qw(t))-(2-it)\right|\leq\frac{C}{N^{1/6}}.
\end{equation}
Hence, for $|t|\leq N^{1/3}$ we have the bound 
\begin{equation}
\label{eq:k10}
\left |
\left(\frac{1}{a\sqrt{N/2}}\right)^{r-s}\frac{(1-qw(t))^s}{(1-qz(u))^r}
\frac{(w(t)-q)^M}{(z(u)-q)^M}\right|\leq C,
\end{equation}
and we see that the part of the integral in (\ref{eq:k4})
where $|t|\leq N^{1/3}$ has a uniform bound for $\xi$, $\eta$
in a compact set.
This proves claim (2) in lemma~\ref{thm:pnglimit} for $r\leq s$.

It also follows from~(\ref{eq:k2}),~(\ref{eq:k3}),~(\ref{eq:k4p}),
~(\ref{eq:k5}),~(\ref{eq:k6}),~(\ref{eq:k9}),~(\ref{eq:k10})
and the dominated convergence theorem that the integral in~(\ref{eq:k4})
converges to 
\begin{equation}
\label{eq:k11}
\frac{\sqrt{2^{\svariable-\rvariable}e^{\eta^2-\xi^2}}}{2\pi^2 i}
\int_{\gamma_1} du 
\int_{\mathbb{R}} dt\,
\frac{1}{2-it-u}
e^{2\xi u-u^2}
e^{(2-it)^2-2\eta(u-it)} \frac{(2-it)^s}{u^r}.
\end{equation}
Now let $v=2-it$. Then we should integrate
$v$ along the line $\Real v=2$ from minus to plus infinity,
call this contour $\Gamma'$.
We obtain the integral
\begin{equation*}
\frac{\sqrt{2^{\svariable-\rvariable}e^{\eta^2-\xi^2}}}{2(\pi i)^2}
\int_{\gamma_1} du \, 
\frac{e^{2\xi u-u^2}}{u^{\rvariable}}
\int_{\Gamma'} dv \,
\frac{v^\svariable}{v-u} e^{v^2-2\eta v}.
\end{equation*}
Expand $(v-u)^{-1}$ as a geometric series. This turns the 
expression into 
\begin{equation}
\label{eq:18}
\frac{\sqrt{2^{\svariable-\rvariable}e^{\eta^2-\xi^2}}}{2(\pi i)^2}
\sum_{k=0}^\infty
\int_{\gamma_1} du \, 
\frac{e^{2\xi u-u^2}}{u^{\rvariable-k}}
\int_{\Gamma'} dv \,
v^{\svariable-k-1} e^{v^2-2\eta v}
\end{equation}
and we recognize the classical  integral representations of the 
Hermite polynomials. The expression now becomes
\begin{equation}
\label{eq:19}
\sum_{k=0}^{\infty}
\sqrt{\frac{(\svariable-k-1)!}{(\rvariable-k-1)!}}
h_{\rvariable-k-1}(\xi)h_{\svariable-k-1}(\eta)\sqrt{e^{-\xi^2-\eta^2}},
\end{equation}
which proves claim (1) in the lemma in the case $r\leq s$.

We now turn our attention to the case $r>s$. Deforming the
$w$-contour through the $z$-contour in~(\ref{eq:4}), 
we get the same integral as above
save for a residue that we pick up at $z=w$. This is
\begin{equation}
\label{eq:10}
\frac{1}{2\pi} 
\int_{-\pi}^\pi 
e^{i(y-x)\theta}
\frac{1}{(1-qe^{i\theta})^{r-s}} \, d\theta.
\end{equation}
We see that the argument above goes through for the remaining
integral also when $r>s$ until~(\ref{eq:18}).
From there the terms $k=0$, \dots, $s+1$ give~(\ref{eq:19})
as before. In terms $k=s+2$, \dots, $r+1$ we instead evaluate
the $v$-integral using the formula 
\begin{equation}
\label{eqn:magic}
\frac{1}{\pi i} \int_{\Gamma} \frac{e^{v^2-2\eta v}}{v^n} \,dv=
\frac{2^n}{\sqrt{\pi}(n-1)!}\int_{\eta}^\infty (\xi-\eta)^{n-1}e^{-\xi^2}\,d\xi,
\end{equation}
which is valid for $n\geq 1$ and will be proved below. 
That accounts for the terms $j=-r$, \dots, $s+1$ in 
definition~\ref{def:herm}. 

Using the well known formula 
\begin{equation*}
\frac{1}{(1-x)^n}
=\sum_{k=0}^\infty
\binom{n+k-1}{k}
x^k
\end{equation*}
the integral in~(\ref{eq:10}) becomes 
\begin{equation*}
\frac{1}{2\pi}
\sum_{k=0}^\infty \binom{r-s+k-1}{k}q^k
\int_{-\pi}^\pi
e^{i(y-x+k)\theta}
\,d\theta.
\end{equation*}
It is readily solved as 
\begin{equation*}
\begin{cases}
\binom{r-s+x-y-1}{x-y}q^{x-y} &\text{if $y\leq x$}\\
0 &\text{if $y> x$.}
\end{cases}
\end{equation*}
With our rescaling, $x=aN+\xi\sqrt{2bN}$
and $y=aN+\eta\sqrt{2bN}$ and the factors $g(r,\xi,N)/g(s,\eta,N)$
we see that the integral in ~(\ref{eq:10}) is
\begin{multline}
\label{eq:9}
\sqrt{2^{s-r}e^{\eta^2-\xi^2}}
q^{(\eta-\xi)\sqrt{2bN}}
(\frac{q}{a\sqrt{N/2}})^{r-s}
\sqrt{2bN}  \times \\\frac{\Gamma(r-s+(\xi-\eta)\sqrt{2bN})}{
\Gamma((\xi-\eta)\sqrt{2bN})
(r-s-1)!} q^{(\xi-\eta)\sqrt{2bN}}
H(\xi-\eta)
\end{multline}
where $H$ is the Heaviside function.
As $N\rightarrow\infty$ we get the limit
\begin{equation}\label{eq:7}
\sqrt{e^{\eta^2-\xi^2}2^{r-s}}
\frac{(\xi-\eta)^{r-s-1}}{(r-s-1)!}
H(\xi-\eta),
\end{equation}
at least for $\xi\neq \eta$. The case $\xi=\eta$ is
a set of measure zero and is not important.
Together with the result for the double integral this 
completes the proof of claim (1). 
It remains to show the estimate in claim (2) in this case. 
But this is easy. 
The expression in  (\ref{eq:9}) is the exact solution 
of  integral~(\ref{eq:10}), and since this is bounded in $N$ for 
$\xi$, $\eta$ in a compact set, 
claim (2) follows.

It remains to show the formula~(\ref{eqn:magic}). 
Observe that, by repeated partial integration,
\begin{equation}
\int_{\eta}^\infty (\xi-\eta)^{n-1} e^{-2\xi v} \,d\xi=
\frac{(n-1)!e^{-2\eta v}}{(2v)^n}
\end{equation}
if $\eta>0$ and $v\in\Gamma'$.
So in this case the left hand side of our formula can be written 
\begin{equation}
\frac{2^n}{\pi i(n-1)!} \int_{\Gamma'} \int_\eta^\infty (\xi-\eta)^{n-1}e^{v^2-2\xi v}
\,d\xi\, dv.
\end{equation}
Here we can change the order of integration and evaluating the Gaussian integral 
gives the right hand side of (\ref{eqn:magic}).

When $\eta<0$, make a change of variables $v\mapsto -v$. The left
hand side of (\ref{eqn:magic}) becomes 
\begin{equation}
\frac{(-1)^n}{\pi i} \int_{\Gamma''} \frac{e^{v^2+2\eta v}}{v^n} \,dv,
\end{equation}
where $\Gamma''$ is parameterised $v=-2+it$, $t=-\infty\mapsto \infty$. 
By deforming the contour $\Gamma''$ into $\Gamma'$ we get
\begin{equation}
\frac{(-1)^{n-1}}{\pi i}\int_\gamma \frac{e^{v^2+2\eta v}}{v^n} \,dv
+\frac{(-1)^n}{\pi i} \int_{\Gamma'} \frac{e^{v^2-2(-\eta) v}}{v^n} \,dv,
\end{equation}
where $\gamma$ is a circle around the origin.
The right term can be evaluated using the the
result for $\eta>0$ above. The left term can be rewritten using the equality between the two
integral formulas for the Hermite polynomials mentioned above and we obtain
\begin{equation}
\frac{2^n}{\sqrt\pi(n-1)!}\int_{-\infty}^\infty 
(\xi-\eta)^{n-1}e^{-\xi^2}\,d\xi-
\frac{2^n}{\sqrt\pi(n-1)!}\int_{-\infty}^\eta
(\xi-\eta)^{n-1}e^{-\xi^2}\,d\xi,
\end{equation}
which proves formula (\ref{eqn:magic}) for $\eta<0$.
\end{proof}

\begin{proof}[Proof of lemma~\ref{krawtchouk_to_hermite}]
Assume first that  $r\leq s$. 
By proposition~\ref{thm:aztec} we have to consider the integral
\begin{equation*}
\frac{\sqrt{N/2}}{(2\pi i)^2}
\int_{\gamma_{r_2}} dz
\int_{\gamma_{r_1}} \frac{dw}{w}
\frac{1}{z-w}
 e^{N(f(z)-f(w))}
\frac{w^{s+\sqrt{N/2}\eta}}{z^{r+\sqrt{N/2}\xi}}
\frac{(1-w)^s}{(1-z)^r}
\frac{(1+z)^{r}}{(1+w)^{s}},
\end{equation*}
where $\gamma_r$ is a circle around $-1$
with radius $r$ oriented anticlockwise, 
$1<r_1<r_2<2$ and
\begin{equation*}
f(z)=\frac{1}{2}\ln z -\ln(1+z).
\end{equation*}
(Here we have ignored the difference between 
$N/2+\xi \sqrt{N/2}$ and its integer part.)
In the proof of lemma~\ref{thm:pnglimit} we could chose the contours of
integration as circles centred at the origin. This cannot be done here.

Note that $f'(z)=0$ gives $z=1$. This leads us to choose
\begin{equation*}
r_1=2-\frac{2}{\sqrt{N/8}}
\end{equation*}
and to deform $\gamma_{r_2}$ to a circle $\Gamma$
oriented clockwise around $1$ with radius $1/\sqrt{N/8}$.
The specific choice of radii are convenient for the computations below.
Choose 
\begin{equation*}
g(r,\xi,N)=\sqrt{N^{-r}e^{-\xi^2}}.
\end{equation*}
Then,
\begin{multline}
\label{eq:l1}
\frac{g(r,\xi,N)}{g(s,\eta,N)}
\sqrt{N/2}K^{\Aztec}(r,\lfloor N/2+\xi\sqrt{N/2}\rfloor;\,
s,\lfloor N/2+\eta\sqrt{N/2}\rfloor)\\
=\sqrt{N^{s-r}e^{\eta^2-\xi^2}}
\frac{\sqrt{N/2}}{(2\pi i)^2}
\int_{\Gamma}dz\int_{\gamma_{r_1}} \frac{dw}{w}
\frac{1}{z-w}
e^{N(f(z)-f(w))}
\frac{w^{s+\eta\sqrt{N/2} }}{z^{r+\xi\sqrt{N/2} }}
\frac{(1-w)^s}{(1-z)^r}
\frac{(1+z)^r}{(1+w)^s}.
\end{multline}
Parameterize $\gamma_{r_1}$ by 
\begin{equation}
\label{eq:l1p}
w(t)=-1+r_1e^{itE_N},
\end{equation}
for $-\pi/E_N\leq t\leq \pi/E_N$,
$E_N=1/\sqrt{N/2}$.
We have
\begin{align*}
\Real(f(w(0))-f(w(t)))&=\frac{1}{2}\ln\left|\frac{w(0)}{w(t)}\right|\\
&=-\frac{1}{4}\ln\left(1+\frac{2r_1(1-\cos E_N t)}{(r_1-1)^2}\right)\\
&\leq-\frac{1}{4} \ln \left(1+\frac{1}{2} (1-\cos E_Nt)\right)
\intertext{ for large enough $N$. 
Again $\cos x\leq 1-x^2/8$ when $|x|\leq \pi$, the last expression is}
&\leq -\frac{1}{4}\ln\left(1+E_N^2t^2/16\right)\leq -Ct^2/N
\end{align*}
for $|t|\leq \pi/E_N$, where $C> 0$ is an absolute constant.
Hence,
\begin{equation}
\label{eq:l2}
\Real (f(w(0))-f(w(t)))\leq -Ct^2/N
\end{equation}
for $|t|\leq \pi/E_N$, with $C>0$.

In the right hand side of~(\ref{eq:l1}) we make the change of variables
\begin{equation}
\label{eq:l3}
z=z(u)=1-u/\sqrt{N/8}
\end{equation}
with $u$ on the unit circle oriented anticlockwise, denoted $\gamma$. 
We obtain the integral 
\begin{multline}
\label{eq:l4}
\sqrt{N^{s-r}e^{\eta^2-\xi^2}}
\frac{iE_N}{2\pi^2}
\int_\gamma du
\int_{-\pi/E_N}^{\pi/E_N} dt
\frac{1}{z(u)-w(t)}
e^{N(f(z(u))-f(w(t)))}\\
\times
\frac{(w(t))^{s-1+\eta\sqrt{N/2} }}{(z(u))^{r+\xi\sqrt{N/2} }}
\frac{(1-w(t))^s}{(1-z(u))^r}
\frac{(1+z(u))^r}{(1+w(t))^{s-1}}.
\end{multline}

Note that 
\begin{equation}
\label{eq:l4p}
f(1+h)=f(1)-h^2/8+O(h^3)
\end{equation}
for small $|h|$. Hence, for $N$ sufficiently large
\begin{equation}
\label{eq:l5}
N(f(z(u))-f(w(0)))=-u^2+4+h_N(u)/\sqrt{N},
\end{equation}
where $h_N(u)$ is bounded for $|u|=1$. We have
\begin{equation}
\label{eq:l6}
\left|\frac{(w(t))^{s-1+\eta\sqrt{N/2}}}{(z(u))^{r+\xi\sqrt{N/2}}}\right|
\leq 3^{s-1+\eta\sqrt{N/2}}\left(1-\frac{u}{\sqrt{N/8}}\right)^{-r-\xi\sqrt{N/2}}\leq C 3^{s+\eta\sqrt{N/2}}
\end{equation}
for some constant $C>0$ depending on $r$.
We also have
\begin{equation}
\label{eq:l7}
\sqrt{N/8}|z(u)-w(t)|\geq 1
\end{equation}
for $u\in \gamma$ and $|t|\leq \pi/E_N$, 
and
\begin{equation}
\label{eq:l8}
\left|\sqrt{N^{s-r}} \frac{(1-w(t))^s}{(1-z(u))^r} \frac{(1+z(u))^r}{(1+w(t))^{s-1}}\right|\leq C N^{s/2}
\end{equation}
for $u\in \gamma$, $|t|\leq \pi/E_N$, by~(\ref{eq:l3}) and~(\ref{eq:l1p}).

It follows from~(\ref{eq:l2}),~(\ref{eq:l5}),~(\ref{eq:l6}),~(\ref{eq:l7})
and~(\ref{eq:l8}) that the part
of the integral~(\ref{eq:l4}) where the $t$-integration is restricted to 
$N^{1/3}\leq |t|\leq \pi/E_N$ can be bounded by
\begin{equation*}
CN^{s/2}3^{s+\eta\sqrt{N/2}}\int_{|t|\geq N^{1/3}} e^{-Ct^2}\,dt,
\end{equation*}
which tends to $0$ as $N\rightarrow\infty$. When
$|t| \leq N^{1/3}$ and $u\in\gamma$, we have 
\begin{equation}
\label{eq:l9}
\left|\frac{(w(t))^{s-1+\eta\sqrt{N/2}}}{(z(u))^{r+\xi\sqrt{N/2}}}\right|
\leq C 
\end{equation}
and
\begin{equation}
\label{eq:l10}
\left|\sqrt{N^{s-r}} \frac{(1-w(t))^s}{(1-z(u))^r} \frac{(1+z(u))^r}{(1+w(t))^{s-1}}\right|\leq C,
\end{equation}
where $C$ depends on $s$, $r$ and $\eta$ but is independent of $N$.

Hence, we see that the part of the integral in~(\ref{eq:l4})
where $|t|\leq N^{1/3} $ has a uniform bound for $\xi$ and $\eta$ in 
a compact set. This proves claim (2) for $r\leq s$. 

It also follows 
from~(\ref{eq:l2}),~(\ref{eq:l3}),~(\ref{eq:l4p}),~(\ref{eq:l5}),\
~(\ref{eq:l7}),~(\ref{eq:l8}),~(\ref{eq:l9}),~(\ref{eq:l10}) and the dominated 
convergence theorem that the integral in (\ref{eq:l4})
converges to 
\begin{equation*}
\sqrt{2^{s-r}e^{\eta^2-\xi^2}}\frac{i}{2\pi^2}
\int_\gamma du\int_{\mathbb{R}} dt
\frac{1}{(2-it)-u}
e^{(2-it)^2-2(2-it)\eta}e^{-u^2+2u\xi}
\frac{(2-it)^s}{u^r},
\end{equation*}
which is exactly the integral in~(\ref{eq:k11}). This proves 
claim~(1) in the lemma in the case $r\leq s$.

For $r>s$ we can deform the contours one through the other
to get the same integral as we solved above. On the way
we pick up the residue of a pole at $z=w$. It is
\begin{equation}
\label{eq:l11}
\frac{1}{2\pi i}\int_{\gamma} 
\frac{dw}{w} w^{-(r-s)-(x-y)} \left(\frac{1+w}{1-w}\right)^{r-s},
\end{equation}
where $x=\lfloor N/2+\xi\sqrt{N/2}\rfloor$ 
and $y=\lfloor N/2+\eta\sqrt{N/2}\rfloor$.
The argument above goes through for the remaining integral
also when $r>s$.
We see that if $\eta>\xi$, then $x-y\rightarrow -\infty$ and this last
 integral is zero.
For simplicity, let $k=r-s$ and $\beta=(x-y)/\sqrt{N/2}$.
The coefficient in front of $w^j$ in the expansion of $[(1+w)/(1-w)]^k$
is 
\begin{equation*}
\frac{1}{2\pi i}\int_{\gamma} \frac{dw}{w} w^{-j}\left(\frac{1+w}{1-w}\right)^k=
\sum_{i=0}^k \binom{k}{i} \binom{j+i-1}{j} (-1)^{k-i} 2^i.
\end{equation*}
One then sees that the $i=k$ term dominates 
when $N$ is large. 
\begin{equation}
\label{eq:l13}
\left|\sum_{i=0}^{k-1} \binom{k}{i} \binom{j+i-1}{j} (-1)^{k-i} 
2^i\right|\leq CN^{k-1}.
\end{equation}
Keeping only the $i=k$ term and 
plugging in  our rescaling and the factors $g(r,\xi,N)/g(s,\eta,N)$,
we see that the integral in (\ref{eq:l11}) is
\begin{equation}
\label{eq:l14}
\sqrt{e^{\eta^2-\xi^2} N^{s-r}}\sqrt{N/2}\,
2^{r-s}\binom{\beta\sqrt{N/2}+2(r-s)-1}{r-s-1}H(\xi-\eta)
\end{equation}
for $\xi\geq \eta$. When $\xi\neq \eta$ this tends to 
\begin{equation*}
\label{eq:l12}
\sqrt{e^{\eta^2-\xi^2}2^{r-s} }\frac{(\xi-\eta)^{r-s-1}}{(r-s-1)!}H(\xi-\eta)
\end{equation*}
as $N\rightarrow\infty$ which together with the corresponding result for 
the double integral settles claim (1) in the case $r>s$.

Claim (2) in this case follows from the corresponding result for
the double integral,(\ref{eq:l13}) and the boundedness of the
expression in~(\ref{eq:l14}).
\end{proof}
\appendix
\section{Asymptotics for Hahn polynomials}
\label{sec:asympt-hahn-polyn}
The Hahn polynomials, as they are defined in~\cite{koekoek:askey},
satisfy 
\begin{equation*}
\sum_{x=0}^N \binom{\alpha+x}{x}\binom{\beta+N-x}{N-x}
Q_m(x; \alpha,\beta,N) Q_n(x; \alpha,\beta,N)
= (d^{(\alpha,\beta)}_{n,N})^2\delta_{nm}
\end{equation*}
where 
\begin{equation*}
(d^{(\alpha,\beta)}_{n,N})^2=
\frac{(-1)^n (n+\alpha+\beta+1)_{N+1}(\beta+1)_n n!}{
(2n+\alpha+\beta+1)(\alpha+1)_n(-N)_nN!}.
\end{equation*}
The Hermite polynomials are defined as usual:
\begin{equation}
\frac{1}{\sqrt{\pi}}\int_{\mathbb{R}} H_n(x) H_m(x)\,dx=2^n n!\delta_{nm}.
\end{equation}
With this notation, the following well known limit theorem holds.
\begin{thm}
Let $0<p<1$ and $\gamma\geq 0$. 
Let $\tilde x= \lfloor pN+x\sqrt{2p(1-p)N(1+\gamma^{-1})}\rfloor $ and
\begin{equation*}
f_{n,N}=(-1)^n\sqrt{\binom{N}{n}2^n n!\left(\frac{p}{1-p}\right)^n
\left(\frac{\gamma}{1+\gamma}\right)^n} 
\end{equation*}
\begin{equation}
\label{eq:11}
E_n(x)=f_{n,N}
Q_n(\tilde x ; \gamma pN, \gamma(1-p)N,N).
\end{equation}
Then 
\[E_n(x)-H_n(x) = O(\sqrt{N^{-1}})\]
uniformly on compact sets.
\end{thm}
\begin{proof}
The idea is induction on $n$.
To start with, $Q_0(y, \alpha,\beta,N)=1$ and we actually 
have $E_0(x)=H_0(x)$. For $n=1$, 
\[
Q_1(y,\alpha,\beta,N)=1-\frac{2+\alpha+\beta}{(\alpha+1)N}x
\]
so 
\begin{align*}
E_1(x)&=-\sqrt{ 2N\left(\frac{p}{1-p}\right)
\left(\frac{\gamma}{1+\gamma}\right)}
\left(1-\frac{2+\gamma N}{(\gamma pN+1)N}
\left(pN+\sqrt{2p(1-p)N(1+\gamma^{-1})} x \right)\right)\\
&=\cdots = H_1(x)+O(\sqrt{N^{-1}}).
\end{align*}
Now assume that the theorem is true for $n$ and $n-1$. We wish to
show that it is true for $n+1$.

There are three term recursion formulas for both 
Hahn and Hermite polynomials.
Let 
\begin{align*}
A_n&=\frac{(n+\alpha +\beta+1)(n+\alpha+1)(N-n)}{
(2n+\alpha+\beta+1)(2n+\alpha+\beta+2)}\\
C_n&=\frac{n(n+\alpha+\beta+N+1)(n+\beta)}{
(2n+\alpha+\beta)(2n+\alpha+\beta+1)}.
\end{align*}
Then
\begin{align}
\label{eq:12}
A_nQ_{n+1}(x) &=(A_n+C_n-x)Q_n(x)-C_nQ_{n-1}(x)\\
H_{n+1}(x)&=2xH_n(x)-2nH_{n-1}(x).
\end{align}
Solving~(\ref{eq:11}) for $Q_n$ and inserting into~(\ref{eq:12})
gives after some simplification the following:
\begin{equation}
\label{eq:17}
E_{n+1}(x) =
\frac{f_{n+1,N}}{f_{n,N}}\left(1+\frac{C_n}{A_n}-\frac{\tilde x}{A_n} \right)
E_n(x)
-\frac{f_{n+1,N}}{f_{n-1,N}}\frac{C_n}{A_n} E_{n-1}(x).
\end{equation}
Observe that under our scaling, 
\begin{align*}
A_n&=pN+O(N^{-1})\\
C_n&=\frac{(1+\gamma)n(1-p)}{\gamma}+O(N^{-1}).\\
\end{align*}
Inserting this into equation~(\ref{eq:17}) and doing some manipulations
gives
\begin{align*}
E_{n+1}(x)& =\left( 2x+O(N^{-1/2})\right) E_n(x)
+\left( 2n+O(N^{-1/2})\right) E_{n-1}(x),
\intertext{which with our induction assumption is}
&= 2xH_n(x)+2nH_{n-1}(x) + O(N^{-1/2})\\
&= H_{n+1}(x)  + O(N^{-1/2}).
\end{align*}
This completes the proof.
\end{proof}
Applying  Stirling's approximation to $d_{n,N}^{\alpha,\beta}$, 
$f_{n,N}$ and the weight function $w^{(\alpha,\beta)}_N(x)$, it is easy 
to show that 
\begin{cor}
As before, $\tilde x=pN+x\sqrt{2p(1-p)N(1+\gamma^{-1})}$.
\begin{equation}
\sqrt[4]{2p(1-p)N(1+\gamma^{-1})}q^{(\alpha,\beta)}_{n,N}(\tilde x)
\sqrt{w^{\alpha,\beta}_{N}(\tilde x)}
\longrightarrow
(-1)^n h_n(x) e^{-x^2/2}
\end{equation}
as $N\rightarrow\infty$ if $\alpha/N\rightarrow p\gamma$ and
$\beta/N\rightarrow (1-p)\gamma$.
\end{cor}
\clearpage
\bibliography{thesisrapport}

\end{document}